\documentclass[a4paper,oneside,11pt]{amsart}

\reversemarginpar

\usepackage[colorlinks,hyperindex,linkcolor=blue,urlcolor=black,pdftitle={Some sufficient conditions for the existence of hyperinvariant subspaces for operators  intertwined with unitaries},pdfauthor={Maria F.Gamal'}]
{hyperref}

\usepackage{amsmath}
\usepackage{amsfonts}
\usepackage{amssymb}
\usepackage{amsthm}

\usepackage[english]{babel}
\usepackage{enumerate}

\usepackage{graphicx}
\usepackage{color}

\usepackage[abbrev]{amsrefs}

\numberwithin{equation}{section}

\theoremstyle{plain}
\newtheorem{theorem}{Theorem}[section]
\newtheorem{lemma}[theorem]{Lemma}
\newtheorem{corollary}[theorem]{Corollary}
\newtheorem{proposition}[theorem]{Proposition}

\theoremstyle{definition}
\newtheorem{definition}[theorem]{Definition}
\newtheorem{remark}[theorem]{Remark}
\newtheorem{example}[theorem]{Example}

\begin{document}

\title[Operators intertwined with unitaries]{Some sufficient conditions for the existence of hyperinvariant subspaces for operators  intertwined with unitaries}

\author{Maria F. Gamal'}
\address{
 St. Petersburg Branch\\ V. A. Steklov Institute 
of Mathematics\\
 Russian Academy of Sciences\\ Fontanka 27, St. Petersburg\\ 
191023, Russia  
}
\email{gamal@pdmi.ras.ru}




\begin{abstract} For a power bounded or  polynomially bounded operator $T$  sufficient conditions for the existence of a nontrivial hyperinvariant subspace are given. The obtained hyperinvariant subspaces of $T$ have the form of the closure of the range of $\varphi(T)$. Here $\varphi$ is a singular inner function, if $T$ is polynomially bounded, or $\varphi$ 
is an analytic in the unit disc function with absolutely summable Taylor coefficients and singular inner part, if $T$ is supposed to  be power bounded only.
 Also, an example of a  quasianalytic contraction $T$ is given such that the quasianalytic spectral set of $T$  is not the whole unit circle $\mathbb T$, while 
$\sigma(T)=\mathbb T$.  Proofs are based on results by Esterle, Kellay, Borichev and Volberg.
\end{abstract}

\maketitle

\section{Introduction}

Let $\mathcal H$ be a (complex, separable) Hilbert space, and let  $\mathcal L(\mathcal H)$ be the algebra of all (linear, bounded) operators acting on  $\mathcal H$. A (closed) subspace  $\mathcal M$ of  $\mathcal H$ is called \emph{invariant} 
for an operator $T$, $T\in\mathcal L(\mathcal H)$, if $T\mathcal M\subset\mathcal M$, and $\mathcal M$ is called \emph{hyperinvariant} 
for $T$ if $C\mathcal M\subset\mathcal M$ for all $C\in\mathcal L(\mathcal H)$ such that $CT=TC$.  The complete lattice of all invariant (resp.,  hyperinvariant) subspaces of $T$ is denoted by  $\operatorname{Lat}T$ (resp., by 
$\operatorname{Hlat}T$). The algebra of all $C\in\mathcal L(\mathcal H)$  such that $TC=CT$ is called the \emph{commutant} of $T$ and is denoted by  $\{T\}'$. 
 The hyperinvariant subspace problem is the question whether for every nontrivial operator 
$T\in\mathcal L(\mathcal H)$ there exists a nontrivial hyperinvariant subspace. 
Here ``nontrivial operator'' means that it is not a scalar multiple of the identity operator, 
and ``nontrivial subspace'' means any subspace different from  $\{0\}$ and  $\mathcal H$.

For Hilbert spaces  $\mathcal H$ and $\mathcal K$, the symbol   $\mathcal L(\mathcal H, \mathcal K)$ means the space of (linear, bounded) 
operators acting from $\mathcal H$ to $\mathcal K$. Suppose that $T\in\mathcal L(\mathcal H)$, $R\in\mathcal L(\mathcal K)$, $X\in\mathcal L(\mathcal H, \mathcal K)$, and $X$  \emph{intertwines} $T$ and $R$, that is, $XT=RX$. If $X$ is unitary, then $T$ and $R$ 
are called  \emph{unitarily equivalent}, in notation: $T\cong R$. If $X$ is invertible, that is, $X^{-1}\in\mathcal L(\mathcal K, \mathcal H)$, 
then $T$ and $R$ are called \emph{similar}, in notation: $T\approx R$.
If $X$ is a \emph{quasiaffinity}, that is, $\ker X=\{0\}$ and $\operatorname{clos}X\mathcal H=\mathcal K$, then
$T$ is called a  \emph{quasiaffine transform} of $R$, in notation: $T\prec R$. If $T\prec R$ and 
$R\prec T$, then $T$ and $R$ are called  \emph{quasisimilar}, in notation: $T\sim R$. 
If $\ker X=\{0\}$, then we write 
 $T \buildrel i \over \prec R$, while
if $\operatorname{clos}X\mathcal H=\mathcal K$, we write $T \buildrel d \over \prec R$. It follows immediately from the definition that 
if $T \buildrel i \over \prec R$ and $\lambda\in\mathbb C$ is an eigenvalue of $T$, then $\lambda$ is an eigenvalue of $R$. Also,  
 $T \buildrel d \over \prec R$ if and only if  $R^\ast \buildrel i \over \prec T^\ast$. Therefore, if $T \buildrel d \over \prec R$ 
and $R^\ast$ has an eigenvalue, then $T^\ast$  has  an  eigenvalue; consequently, $T$ has a nontrivial hyperinvariant subspace.  Note that the relationship $\buildrel d \over \prec$ takes place between 
a power bounded operator and its isometric asymptote (\cite{ker1}; recalled in  Sec. 3 of the present paper). 
Recall that if $T\sim R$ and one of $T$ or $R$ has a  nontrivial hyperinvariant subspace, then another  also  has 
a nontrivial hyperinvariant subspace, see, for example, {\cite[Lemma 2.1]{berc}} or {\cite[Theorem 6.19]{rara}} (and {\cite[Proposition II.5.1]{sznagy}} for a related result).

An  operator $T\in\mathcal L(\mathcal H)$ is called  \emph{power bounded}, if 
$\sup_{n\geq 0}\|T^n\| < \infty$. It is easy to see that for such operators the space 
$$\mathcal H_{T,0}=\{x\in\mathcal H :\ \|T^nx\|\to 0\}$$ is hyperinvariant for $T$ (sf. {\cite[Theorem 6.21]{rara}}, {\cite[Theorem II.5.4]{sznagy}}). Classes $C_{ab}$, where indices $a$ and $b$ can be equal to $0$, $1$, or a dot,  of power bounded operators 
are defined as follows. If $\mathcal H_{T,0}=\mathcal H$, then  $T$ is  \emph{of class} $C_{0\cdot}$, while if  $\mathcal H_{T,0}=\{0\}$, then $T$ is 
 \emph{of class} $C_{1\cdot}$. Furthermore,  $T$  is \emph{of class} $C_{\cdot a}$, if $T^\ast$ is of class  $C_{a\cdot}$,  
 and $T$ is  \emph{of class} $C_{ab}$, if $T$ is of classes $C_{a\cdot}$ and $C_{\cdot b}$, $a$, $b=0,1$. 

 The operator $T\in\mathcal L(\mathcal H)$ is called  \emph{polynomially bounded}, if there exists a constant $M$ such that 
$\|p(T)\|\leq M \sup\{|p(z)|:  |z|\leq 1\} $ for every polynomial $p$. For a polynomially bounded operator $T\in\mathcal L(\mathcal H)$ 
there exist $\mathcal H_a$, $\mathcal H_s\in\operatorname{Hlat}T$ such that $ \mathcal H=\mathcal H_a\dotplus\mathcal H_s$,   $T|_{\mathcal H_a}$ is an 
 \emph{absolutely continuous (a.c.)} polynomially bounded operator, and $T|_{\mathcal H_s}$ is similar to a 
singular unitary operator. Thus, if $\mathcal H_s\neq\{0\}$, then $T$ has nontrivial  hyperinvariant subspaces. The  definition of a.c. polynomially bounded operators is not recalled here,  because it will be not used in the present paper. We recall only that \emph{$T$ is an  a.c. polynomially bounded operator if and only if $T$ admits an $H^\infty$-functional calculus}
  \cite{mlak}, {\cite[Theorem 23]{ker5}}. (Although many results on polynomially 
bounded operators that will be used in the present paper were originally  proved by Mlak, we will refer to \cite{ker5} for the convenience of references.)
For the existence of invariant subspaces of polynomially 
bounded operators see \cite{rej}.

The operator $T\in\mathcal L(\mathcal H)$  is called 
a  \emph{contraction}, if $\|T\|\leq 1$. 
A contraction is polynomially bounded with the constant 
$1$ (von Neumann inequality; 
see, for example, {\cite[Proposition I.8.3]{sznagy}} or {\cite[Proposition X.1.7]{beau}}). Clearly, a polynomially bounded operator is power bounded. 
(It is well known that the converse is not true, see \cite{fog} for the first example of power bounded and not polynomially bounded operator, 
and \cite{pis} for  the first example of polynomially bounded operator which is not similar to a contraction.)

Although the hyperinvariant subspace problem is reduced to the case of some subclass of contractions of class $C_{00}$ \cite{foha}, the case 
of power bounded operators that are not of class $C_{00}$ seems to be more tractable. For example, if $T$ is a power bounded operator of class $C_{11}$, 
then $T$ is quasisimilar to a unitary operator and, consequently,  has  nontrivial  hyperinvariant subspaces {\cite[Theorem 6.20]{rara}},
{\cite[Proposition II.5.3]{sznagy}} (see {\cite[Ch.  IX.4]{sznagy}} and \cite{ker1} for the further studies).

\subsection{Main results}

In the present paper, we consider power bounded operators  that are quasiaffine transforms of a.c. unitaries. Some sufficient conditions for the existence of nontrivial hyperinvariant subspaces of such operators are given. For $T\in\mathcal L(\mathcal H)$ these   hyperinvariant subspaces are the closures of $\operatorname{ran}\varphi(T):=\varphi(T)\mathcal H$. The operator $\varphi(T)$ is the function of $T$ obtained using appropriate functional calculus. For absolutely continuous (a.c.) polynomially bounded operators,  $H^\infty$-functional calculus is used. For power bounded operators,  functional calculus on the algebra 
$A^+(\mathbb T)$ of functions analytic in the unit disc with absolutely summable Taylor coefficients is used.  The relation   
$T \buildrel d \over \prec U$, where $U$ is an a.c. unitary operator, is used in the construction. 
The results, in a simple framework, are the following. 

\emph{Suppose that  $T\in\mathcal L(\mathcal H)$,  $U\in\mathcal L(\mathcal K)$, 
 $X\in\mathcal L(\mathcal H,\mathcal K)$, $U$ is an a.c. unitary operator, $XT=UX$, and 
$\operatorname{clos}X\mathcal H=\mathcal K$. If there exists $0\neq g\in\mathcal K$ such that $\|T^{\ast n}X^\ast g\|$ tends to $0$ sufficiently fast, then there exists a function $\varphi$ such that 
$\operatorname{clos}\varphi(T)\mathcal H\neq\mathcal H$.}

 Note that we can suppose that  $\|T^{\ast n} x\|\to 0$ for every $x\in\mathcal H$. (Otherwise, the existence of a  nontrivial hyperinvariant subspace for $T$ follows from {\cite[Theorem 6.21]{rara}} or {\cite[Theorem II.5.4]{sznagy}}.)

Our main results are Theorems \ref{T:6.7} and \ref{T:6.9} (for polynomially bounded and power bounded operators, respectively). Corollaries  \ref{C:6.8} and \ref{C:6.10} give simpler sufficient conditions for operators to satisfy the conditions of Theorems \ref{T:6.7} and \ref{T:6.9}.

Our proofs  are based on the results in \cite{est} and \cite{kell1}.

Furthermore,  an example of a quasianalytic contraction $T$ such that $\sigma(T)=\mathbb T$ and  $\pi(T)\neq\mathbb T$, where 
 $\pi(T)$ is the quasianalytic spectral set of $T$, is given (Theorem~\ref{T:7.7}).
(See \cite{ker5} and references therein for the definition and discussions; see also Sec. 1.3 of the present paper.) 
The relationship  $\pi(T)\subset\sigma(T)\cap\mathbb T$ is always true. Examples of contractions $T$ such that  
$\pi(T) = \sigma(T)\cap\mathbb T\neq\mathbb T$ are known (see \cite{kers}
 and references therein and in \cite{ker5}). However, in these examples $\sigma(T)\not\subset\mathbb T$, and, as far as
the author understand, the known methods do not allow to construct examples of $T$ such that $\emptyset\neq\pi(T)\neq\mathbb T$ and $\sigma(T)\subset\mathbb T$. 
On the other hand, the method from the present paper does not allow to construct  examples of $T$ such that 
$\mathbb T\not\subset\sigma(T)$. Thus,  the question  from \cite{kers} whether there exists a
 contraction $T$ such that $\pi(T)=\sigma(T)\neq\mathbb T$ remains open. Recall that  contractions $T$ such that $\pi(T)=\sigma(T)=\mathbb T$ can be found among weighted shifts, see \cite{est} and  \cite{kers}.  The proof of the quasianalyticity of $T$ constructed in the present paper  is based on \cite{bovo}. 
 The author does not known, whether this 
 $T$ satisfies the conditions of  Theorem~\ref{T:6.7} of the present paper. Therefore, to find  nontrivial hyperinvariant subspaces of $T$, the results from \cite{kell3} are used.

\subsection{Notations and auxiliary results}

Symbols $\mathbb D$, $\operatorname{clos}\mathbb D$, and $\mathbb T$ denote the open unit disc, the closed unit disc, 
and the unit circle, respectively. The normalized Lebesgue measure on $\mathbb T$ is denoted by $m$. Furthermore,  $L^2=L^2(\mathbb T,m)$, symbol $\chi$ denotes the identity function, i.e., $\chi(\zeta)=\zeta$, 
$\zeta\in\mathbb T$.  Clearly, $\{\chi^n\}_{n\in\mathbb Z}$ is an orthonormal basis of $L^2$. Symbol  $U_{\mathbb T}$  denotes  the operator of multiplication by $\chi$ acting on $L^2$.
 For a Borel set $\tau\subset\mathbb T$, the restriction of $U_{\mathbb T}$ on its reducing subspace $L^2(\tau,m)$ will be denoted  by $U(\tau)$. For $1\leq p\leq\infty$, $H^p$ is the Hardy space on $\mathbb D$. The space $H^2$ can be regarded as a subspace of $L^2$, put $H^2_-=L^2\ominus H^2$. 
Symbols $P_+$ and $P_-$ denote the orthogonal projections  on $H^2$ and $H^2_-$, respectively. 
Symbols $S$ and $S_\ast$ denote the restriction and the compression of $U_{\mathbb T}$ on $H^2$ and $H^2_-$, respectively. 
Clearly, $U_{\mathbb T}$ has the following form with respect to the decomposition $L^2 = H^2\oplus H^2_-$:
$$U_{\mathbb T}=\begin{pmatrix} S & (\cdot, \chi^{-1})\chi^0 \\ \mathbb O & S_\ast\end{pmatrix}.$$

For a function $\varphi$ analytic in $\mathbb D$ set  $\widetilde{\varphi}(z)=\overline{\varphi(\overline z)}$, $z\in\mathbb D$. Clearly, 
$\widetilde{\varphi}$ is analytic in $\mathbb D$, and $\widehat{\widetilde\varphi}(n)=\overline{\widehat{\varphi}(n)}$, $n\geq 0$. If $\varphi\in H^\infty$ and  $T$ is  an a.c. polynomially bounded operator, then  
$\varphi(T^\ast)=\widetilde\varphi(T)^\ast$ (\cite{mlak}, {\cite[Proposition 14]{ker5}}). 

For a   singular inner function $\theta$ let $\mu_\theta$ be a positive finite Borel singular measure 
on $\mathbb T$ such that 
 $$\theta(z)=\exp\int_{\mathbb T}\frac{z+\zeta}{z-\zeta}\text{d}\mu_\theta(\zeta), \ \  \ z\in\mathbb D,$$
 $\operatorname{supp}\mu_\theta$ is the closed support of $\mu_\theta$. 

Symbol $A^+(\mathbb T)$ denotes the Banach algebra of functions analytic in  $\mathbb D$ with absolutely summable Taylor coefficients:
$$A^+(\mathbb T)=\{\varphi\colon\operatorname{clos}\mathbb D\to\mathbb C, \ 
 \varphi(z)\!=\!\sum_{n=0}^\infty\widehat\varphi(n)z^n, 
\ z\in\operatorname{clos}\mathbb D, \ \sum_{n=0}^\infty|\widehat\varphi(n)|<\infty\}.$$
Clearly, $A^+(\mathbb T)\subset H^\infty$.

For a Hilbert space  $\mathcal H$ and a (closed) subspace  $\mathcal M$ of $\mathcal H$,  symbols  $P_{\mathcal M}$  and  $I_{\mathcal H}$ denote the orthogonal projection on $\mathcal M$ and the identity operator on $\mathcal H$, respectively.

Every function $\varphi\in H^\infty$ can be represented as a product of three functions: $\varphi=\varphi_1\varphi_2\varphi_3$, where 
$\varphi_1$ is a Blaschke product, $\varphi_2$ is a singular inner function, and $\varphi_3$ is an outer function 
(some of factors can be absent). It is easy to see 
that if $T$ is an a.c. polynomially bounded operator, and $T^\ast$ has no eigenvalues, then 
$\operatorname{ran} \varphi(T)$ is  dense for every $\varphi\in H^\infty$ without singular inner factor 
(Lemmas~\ref{L:1.2} and \ref{L:1.3} below). Therefore, if for such $T$ there
exists $\varphi\in H^\infty$ such that $\operatorname{ran} \varphi(T)$ is not dense, then there
exists a singular inner function $\theta\in H^\infty$ such that $\operatorname{ran} \theta(T)$ is not dense. 

The following lemma is formulated in \cite{mlak}. It mentioned there that it can proved almost  as 
{\cite[Proposition III.3.1]{sznagy}}  for a.c. contractions. (Of course, the reference in \cite{mlak} given on the first edition of \cite{sznagy}.) For convenience, we give the proof here.

\begin{lemma}[\cite{mlak}]\label{L:1.2}  Suppose that $T$ is an a.c. polynomially bounded operator, and  $\varphi\in H^\infty$ is an outer function.
Then $\ker \varphi(T)=\{0\}$ and  $\operatorname{ran} \varphi(T)$ is  dense.\end{lemma}
\begin{proof} Since $T^\ast$ is also an a.c. polynomially bounded operator (\cite{mlak}, {\cite[Proposition 14]{ker5}}), and  
$\varphi(T)^\ast=\widetilde\varphi(T^\ast)$ for any function $\varphi\in H^\infty$, it is sufficient to prove that 
$\ker \varphi(T)=\{0\}$. 
Denote by $\mathcal H$ the space on which $T$ acts. Let $x\in\mathcal H$. Since $T$ is a.c. polynomially bounded, there exists $f\in L^1(\mathbb T,m)$ 
such that $(\psi(T)x,x)=\int_{\mathbb T}\psi f\text{{\rm d}}m$ for every $\psi\in H^\infty$ (\cite{mlak}, {\cite[Sec.3]{ker5}}). If $\varphi(T)x=0$, then 
$$(T^n\varphi(T)x,x)=\int_{\mathbb T}\chi^n\varphi f\text{{\rm d}}m=0 \ \ \text{ for all } \  n\geq 0. $$ 
 It means that $\varphi f\in H^1$ and $(\varphi f)(0)=0$. Since $\varphi$ is outer, $f\in H^1$ and $f(0)=0$. Therefore, 
$$\|x\|^2=(x,x)= \int_{\mathbb T} f\text{{\rm d}}m =0. $$ \end{proof}

\begin{lemma}\label{L:1.3} Suppose that  $T$ is an a.c. polynomially bounded operator,   $\varphi\in H^\infty$ is a Blaschke product, 
and  $\operatorname{ran} \varphi(T)$ is not dense. Then there exists $\lambda\in\mathbb D$ such that $\varphi(\lambda)=0$ and 
$\overline\lambda$ is an eigenvalue of $T^\ast$.\end{lemma}
\begin{proof}  Set $T_0=T^\ast|_{\ker\widetilde\varphi(T^\ast)}$. Then $\widetilde\varphi(T_0)=\mathbb O$. By \cite{bercpr}, there exists 
an a.c. contraction $R$ such that $T_0\sim R$. It follows that $\widetilde\varphi(R)=\mathbb O$. 
That is, $R$ is a $C_0$-contraction. By {\cite[Proposition III.4.4]{sznagy}}, the minimal function $\psi$ of $R$ divides $\widetilde\varphi$. 
Therefore, $\psi$ is a Blaschke product. By {\cite[Theorem III.5.1]{sznagy}}, if $\lambda\in\mathbb D$ and $\psi(\overline\lambda) =0$, then $\overline\lambda$ is 
an eigenvalue of $R$. Since $T_0\sim R$, we have that $\overline\lambda$ is an eigenvalue of $T_0$, and consequently, of $T^\ast$. Finally, $\varphi(\lambda)=0$, because  $\psi(\overline\lambda) =0$. \end{proof}

The following lemma is very simple, but useful.

\begin{lemma}\label{L:1.1} Suppose  that $T$ and $R$ are a.c. polynomially bounded operators,
$T \buildrel d \over \prec R$ and there exists $\varphi\in H^\infty$ such that $\operatorname{ran} \varphi(R)$ is not dense. 
Then $\operatorname{ran} \varphi(T)$ is not  dense.
\end{lemma}

\begin{proof} Denote by $\mathcal H$ and  $\mathcal K$ the spaces on which $T$  and $R$ act, and  by $X\in\mathcal L(\mathcal H,\mathcal K)$ 
an operator which realizes the relation $T \buildrel d \over \prec R$. We have $XT=RX$. Therefore,   
$ X\varphi(T) = \varphi(R)X$. The latter equality follows from the facts that  polynomials are sequentially dense in 
the weak-$\ast$ topology of $H^\infty$, and the $H^\infty$-functional calculus for $T$ (resp.,  for $R$)  is continuous in the weak-$\ast$ topologies of $H^\infty$ and 
$\mathcal L(\mathcal H)$ (resp.,   $\mathcal L(\mathcal K)$). If $\operatorname{clos}\varphi(T)\mathcal H = \mathcal H$, then 
\begin{align*}\operatorname{clos}\varphi(R)\mathcal K & = 
\operatorname{clos}\varphi(R)\operatorname{clos}X\mathcal H =
\operatorname{clos}\varphi(R)X\mathcal H = \operatorname{clos}X\varphi(T)\mathcal H \\
& =\operatorname{clos}X\operatorname{clos}\varphi(T)\mathcal H =
\operatorname{clos}X\mathcal H=\mathcal K,\end{align*}
a contradiction. \end{proof}

\subsection{Unitary asymptote}

 In \cite{ker5}, the following definition of a unitary asymptote of $T\in\mathcal L(\mathcal H)$ is given: 
a pair $(X,U)$  is called a {\it unitary asymptote} of $T$, if
$U$ is a unitary operator,  $XT=UX$, and for any other pair $(X',U')$
such that $U'$ is a unitary operator and  $X'T=U'X'$ there exists a unique operator $Z$ such that 
$ZU=U'Z$ and $X'=ZX$.  By \cite{ker1} or  {\cite[Sec. 2]{ker5}}, every power bounded operator $T$ has a unitary asymptote.  Moreover,  for power bounded operators this notion coincides with the notion of unitary asymptote 
introduced in \cite{ker1} and constructed using Banach limit   (see Sec. 3 of the present paper for the recalling of this construction). 

Let an operator $T$ have a unitary asymptote $(X,U)$.
As is mentioned in \cite{ker5}, it follows immediately from the definition of a unitary asymptote  that 
 there exists the mapping
$$\gamma_T\colon \{T\}'\to\{U\}', \ \ \gamma_T(C)=D, $$
 where $ D\in\{U\}'$ is  a unique operator  such that $XC=DX$, and $\gamma_T$ is a unital algebra-homomorphism.

The following theorem is an immediate consequence of this fact (sf. {\cite[Lemma IX.1.4]{sznagy}} and \cite{ker1}). 

\begin{theorem}\label{T:3.1} Suppose   that $T\in\mathcal L(\mathcal H)$, and $(X,U)$ is 
 a unitary asymptote of $T$. Then $\sigma(\gamma_T(C))\subset\sigma(C)$ for every $ C\in\{T\}'$, in particular, 
$\sigma(U)\subset\sigma(T)$. Furthermore,  
$$X^{-1}\mathcal N:=\{x\in\mathcal H  \ :  \ Xx\in\mathcal N\} \in\operatorname{Hlat}T$$ for every 
 $\mathcal N\in\operatorname{Hlat}U$. In particular, $\ker X\in\operatorname{Hlat}T$.\end{theorem}

\begin{proof} The statement about  $\{T\}'$ can be proved exactly as in  {\cite[Lemma IX.1.4]{sznagy}} and \cite{ker1}.  
Therefore, its proof is omitted. 

Let $\mathcal N\in\operatorname{Hlat}U$. Clearly, $X^{-1}\mathcal N$ is linear and closed. Suppose that  $C\in\{T\}'$, $D=\gamma_T(C)$, 
and  $x\in X^{-1}\mathcal N$. Then $XCx=DXx\in\mathcal N$. It means that $Cx\in\mathcal N$. Thus,  $X^{-1}\mathcal N\in\operatorname{Hlat}T$.
Since $\ker X = X^{-1}\{0\}$, we conclude that  $\ker X\in\operatorname{Hlat}T$.  \end{proof}

The notion of the quasianalytic spectral set was introduced in \cite{ker5} only for a.c. polynomially bounded operators, while 
this notion can be introduced for every operator  which has  a   unitary asymptote 
 exactly as for a.c. polynomially bounded operators. We do not give the detailed definition of the quasianalytic set 
$\pi(T)$ of an operator $T$ having a unitary asymptote here. The interested readers can do it themselves (or consult with
\cite{kerntuple}). We mention only that if an operator $T\in\mathcal L(\mathcal H)$ has a (nonzero) unitary asymptote $(X,U)$, then $T$ is quasianalytic if and only if  $X^{-1}\mathcal N$ is trivial (i.e., equals to $\{0\}$ or $\mathcal H$) for every 
$\mathcal N\in\operatorname{Hlat}U$. Consequently, if $T$ is not quasianalytic, then $T$ has a nontrivial hyperinvariant subspaces. Furthermore, if a unitary operator $U$ from  the  unitary asymptote $(X,U)$ is a.c., then the residual set 
$\omega(T)$  can be defined exactly as for a.c. polynomially bounded operators. Namely, $\omega(T)$ is turned out the Borel subset of $\mathbb T$ on which the spectral measure of  $U$ is concentrated. The relationship 
$$\pi(T)\subset\omega(T)\subset\sigma(T)\cap\mathbb T$$
follows from the definitions of $\pi(T)$, $\omega(T)$ and Theorem~\ref{T:3.1}. Moreover, $T$ is quasianalytic if and only if 
$\pi(T)=\omega(T)$.

Note that for a.c. polynomially bounded operators there is the relationship 
between quasianalytic and residual sets on the one side and $H^\infty$-functional calculus on the other side 
(see {\cite[Sec. 6]{ker5}}).

\bigskip

The paper is organized as follows. In Sec. 2  we collect some  facts on weighted  shifts to be used later on. 
 In Sec. 3 and 4 sufficient conditions on an operator $T$ for the existence of  singular inner functions $\varphi$ such that $\operatorname{ran} \varphi(T)$  is not  dense are given. 
In Sec. 3  we consider operators with a part similar to a simple unilateral shift, while in Sec. 4  we consider operators intertwined with (may be reductive) a.c. unitaries.  In Sec. 5  an example of a quasianalytic contraction is given. 

\section{Preliminaries: weighted shifts}

In this section we collect some properties of weighted shifts needed for the remaining part of the paper. 
The weighted shifts have been  intensively studied and  these properties are well known to the specialists.

Let $v\colon{\mathbb Z}_+ \to (0,\infty)$ be a nonincreasing function.
Set $$\ell^2_{v+}=\ell^2_{v+}({\mathbb Z_+}) =\big \{u=\{u(n)\}_{n\in{\mathbb Z_+}}:\ \ 
\|u\|_v^2=\sum_{n\in{\mathbb Z_+}}|u(n)|^2v(n)^2<\infty\big\}.$$
The {\it unilateral weighted shift} $S_{v+}\in\mathcal L(\ell^2_{v+})$ acts according to the formula 
$$(S_{v+} u)(n)=u(n-1), \ \ n\geq 1, \ \ \  (S_{v+} u)(0)=0, \ \ u\in\ell^2_{v+}.$$ 
Since $v$ is nonincreasing, $\|S_{v+}\|\leq 1$, that is, 
$S_{v+}$ is a contraction, and it is easy to see that the contraction $S_{v+}$ 
is completely nonunitary, and, consequently, a.c..
If $v(n)=1$ for all $n\in\mathbb Z_+$, we write $\ell^2_{v+}=\ell^2_+$.

Let $\omega\colon\mathbb Z \to (0,\infty)$ be a nonincreasing function.
Set \begin{align*} \ell^2_\omega=\ell^2_\omega(\mathbb Z) & =\big \{u=\{u(n)\}_{n\in\mathbb Z}:\ \ 
\|u\|_\omega^2=\sum_{n\in\mathbb Z}|u(n)|^2\omega(n)^2<\infty\big\}, \\
\ell^2_{\omega +}&=\big \{u\in\ell^2_\omega : \ \ u(n)=0 \ \text{ for } n\leq -1\},\\ 
\ell^2_{\omega -}&=\big \{u\in\ell^2_\omega : \ \ u(n)=0 \ \text{ for } n\geq 0\}.\end{align*}
The {\it bilateral weighted shift} $S_\omega\in\mathcal L(\ell^2_\omega)$ acts according to the formula 
$$(S_\omega u)(n)=u(n-1), \ \ n\in\mathbb Z,  \ \ u\in\ell^2_\omega.$$ 
Clearly, $\ell^2_{\omega +}$ is an invariant subspace of $S_\omega$, and the restriction $S_{\omega +}$ of 
$S_\omega$ on $\ell^2_{\omega +}$ is a unilateral shift.
Also, $\ell^2_{\omega -}$ is a coinvariant subspace of $S_\omega$, and the compression $S_{\omega -}$ of 
$S_\omega$ on $\ell^2_{\omega -}$ acts according to the formula
$$ (S_{\omega -}u)(n) = u(n-1), \ \ n\leq -1,  \ \ u\in\ell^2_{\omega -}. $$
It is easy to see that if $\omega$ and $w$ are two nonincreasing function such that $\omega\asymp w$, then 
$S_\omega\approx S_w$. Also, if there exists $C>0$ such that  $\omega\leq Cw$, then $S_w\prec S_\omega$.

Since $\omega$ is  nonincreasing, $\|S_\omega\|\leq 1$, that is, 
$S_\omega$ is a contraction. It is easy to see that if $\omega$ is  a nonconstant function, then the contraction $S_\omega$ 
is completely nonunitary, and, consequently, a.c.. Also it is easy to see that if  
\begin{align}\label{E:2.1}\sup_{n\in\mathbb Z}\frac{\omega(n-1)}{\omega(n)}<\infty,\end{align}
then  $S_\omega$ is invertible, that is, $S_\omega^{-1}$ is bounded.
The function $\omega$ is called {\it submultiplicative}, if $\omega(n+k)\leq\omega(n)\omega(k)$ for all $n$, $k\in\mathbb Z$.

The function $\omega\colon\mathbb Z \to (0,\infty)$ is called a {\it dissymmetric weight}, if it is nonincreasing, unbounded, 
$\omega(n)=1$ for all $n\geq 0$, $\omega(-n)^{1/n}\to 1$ when $n\to +\infty$, and  $\omega$ satisfies \eqref{E:2.1}. 
For a dissymmetric weight $\omega$, $\ell^2_{\omega +}=\ell^2_+$, and 
 $\sum_{n\in\mathbb Z}|u(n)|^2<\infty$ for every $u\in\ell^2_\omega$, therefore, 
the natural imbedding 
\begin{align}\label{E:2.2} u\mapsto f, \ \ \widehat f(n)=u(n), \ \ \ n\in\mathbb Z, \ \ \ \ \ell^2_\omega\to L^2(\mathbb T,m) \end{align}
is bounded. A dissymmetric weight $\omega$ is called {\it log-concave}, if the sequence 
$$\Bigl\{\frac{\omega(-n-1)}{\omega(-n)}\Bigr\}_{n\geq 0}$$
is nonincreasing. 
If $\omega$ is a log-concave dissymmetric weight, then $\omega$ is submultiplicative {\cite[Proposition 2.4]{est}}. 
If $\omega$ is a submultiplicative dissymmetric weight,
then $\sigma(S_\omega)=\mathbb T$ {\cite[Proposition 2.3]{est}}.   

\begin{theorem}[{\cite[Theorem 5.7]{est}}]\label{T:2.1}  Suppose  that $\omega$ is a dissymmetric weight.
Then there exists a singular inner function $\theta$ such that 
$m(\operatorname{supp}\mu_\theta)=0$ and 
\begin{align}\label{E:2.3} \sum_{n=0}^\infty \frac{1}{\omega(-1-n)^2}\Bigl|\widehat{\frac{1}{\theta}}(n)\Bigl|^2<\infty. \end{align}
Furthermore, $\operatorname{ran}\theta(S_\omega)$ is not dense.\end{theorem}

\begin{lemma}\label{L:2.2}
\begin{enumerate}[\upshape (i)]
 
\item Suppose  that $w$ is a dissymmetric weight, and $\{N_j\}_{j=1}^\infty$ is a sequence of positive integers 
such that $N_1=1$ and $N_j<N_{j+1}$ for all $j\geq 1$. Put $\omega(n)=1$, $n\geq 0$, and $\omega(n)=w_j$, $-N_{j+1}+1\leq n\leq -N_j$, 
$j\geq 1$. Then $\omega$ is a dissymmetric weight.

\item Suppose  that $\{\beta_n\}_{n=0}^\infty$ is a sequence of positive numbers such that $\beta_n\to\infty$. Then there exists 
a dissymmetric weight $\omega$ such that $\omega(-n-1)\leq\beta_n$ for sufficiently large $n$. 

\item Suppose  that $\{\varepsilon_n\}_{n=0}^\infty$ is a sequence of positive numbers such that 
$$\sum_{n=0}^\infty \varepsilon_n^2 < \infty.$$ Then there exists 
a dissymmetric weight $\omega$ such that $$\sum_{n=0}^\infty \varepsilon_n^2\omega(-n-1)^2 <\infty.$$
\end{enumerate}\end{lemma}
\begin{proof} The part (i) of the lemma can be checking straightforward. To prove (ii), set 
$\beta_n '=\inf_{k\geq n-1}\beta_k$ for 
   $n\geq 1$, then $\beta_n '\geq 1$ for sufficiently large $n$. Take an arbitrary dissymmetric weight $w$. There exists a sequence  
  $\{N_j\}_{j=1}^\infty$ of positive integers satisfying (i) and such that $w(-j)\leq \beta_{N_j} ' $ for $j\geq 2$. Let $\omega$ 
be constructed by $w$ and $\{N_j\}_{j=1}^\infty$ as in (i). It is easy to see that $\omega$ satisfies the conclusion of (ii).
To prove (iii), take an arbitrary dissymmetric weight $w$. There exists a sequence  
  $\{N_j\}_{j=1}^\infty$ of positive integers satisfying (i) and such that 
$\sum_{j=1}^\infty w^2(-j)\sum_{n\geq N_j}\varepsilon_n^2 < \infty$. Let $\omega$ 
be constructed by $w$ and $\{N_j\}_{j=1}^\infty$ as in (i). It is easy to see that $\omega$ satisfies the conclusion of (iii).
 \end{proof}

Recall the following definition.
\begin{definition}\label{D:2.3} A closed subset $F$ of $\mathbb T$ is called a {\it Carleson set}, if 
$$\sum_j m(\mathcal I_j)\log m(\mathcal I_j)>-\infty,$$
where $\{\mathcal I_j\}_j$ is the family of disjoint open arcs such that $\mathbb T\setminus F=\cup_j\mathcal I_j$.
\end{definition}

The following theorem is proved in \cite{kell1}, see the proof of $(1)\Rightarrow(2)$ in {\cite[Theorem 2.2]{kell1}}.

\begin{theorem}[\cite{kell1}]\label{T:2.4}  Suppose that  $\{w_n\}_{n\geq 1}$ is a sequence of positive numbers such that the sequences 
$\{w_{n+1}/w_n\}_{n\geq 1}$ and $\{(\log w_n)/n^b\}_{n\geq 1}$ are nonincreasing, $\liminf_nw_n/n^c>0$ for some $b<1/2$ and $c>0$, 
and $$\sum_{n\geq 1}\frac{1}{n\log w_n}<\infty.$$ 
Then there exist a singular inner function $\theta$ and a constant $C>0$ such that $\operatorname{supp}\mu_\theta$ is a Carleson set, 
$m(\operatorname{supp}\mu_\theta)=0$,
 and
$$ \Bigl|\widehat{\frac{1}{\theta}}(n-1)\Bigr|\leq C w_n \ \ \text{ for all } n\geq 1.$$\end{theorem}

\begin{lemma}\label{L:2.5} Suppose  that $\{w_n\}_{n\geq 1}$ is a sequence of positive numbers satisfying the conditions of  Theorem~\ref{T:2.4}.  
\begin{enumerate}[\upshape (i)]
 
\item Put $\omega(n)=1$ and $\omega(-n-1)=w_{n+1}$ for $n\geq 0$. Then $\omega$ is a dissymmetric weight.

\item For every $s>0$ the sequence $\{w_n^s\}_{n\geq 1}$ satisfies the conditions of  Theorem~\ref{T:2.4}, and 
$$\sum_{n\geq 1}\frac{1}{ w_n^s}<\infty.$$ 
\end{enumerate}\end{lemma}
\begin{proof} The part (i)  and the first statement of the part (ii) of the lemma can be checked directly. The second statement of the part (ii) is a consequence of the following statement. For every $\varepsilon>0$ there exists $C>0$ such that $n^{\frac{1}{\varepsilon}}\leq Cw_n$ for sufficiently large $n$. The proof of this statement is contained in the proof of  {\cite[Theorem 3.1]{kell1}}.\end{proof}

\begin{example} [{\cite[Remark 3.2]{kell1}}] For $a>1$ put $w_n=n^{(\log\log n)^a}$ for sufficiently large $n$. Then it is possible to define $w_n$ 
for small $n\geq 1$ in such a way that $\{w_n\}_{n\geq 1}$ satisfies the conditions of Theorem~\ref{T:2.4}.\end{example}

\section{Operators with a part similar to the simple unilateral shift}

For a power bounded operator $T\in \mathcal L(\mathcal H)$ the  {\it isometric asymptote} 
$$\bigl(X_{T, +} \in\mathcal L(\mathcal H,\mathcal H^{(a)}_+),\ T^{(a)}_+\in\mathcal L(\mathcal H^{(a)}_+)\bigr)$$ is defined using Banach limit, 
see \cite{ker1}.  For  convenience, the construction is recalled here.  
 
Let $(\cdot,\cdot)$ be the inner 
product on the Hilbert space $\mathcal H$, 
and let $T\in \mathcal L(\mathcal H)$  be a  power bounded operator.
 Define a new semi-inner product on
$\mathcal H$ by the formula $\langle x,y\rangle=\operatorname{Lim}_{n\to\infty}(T^nx,T^ny)$,
 where $x,y\in \mathcal H$, and $\operatorname{Lim}$ is  some fixed Banach limit. Set 
$\mathcal H_0=\mathcal H_{T,0}=\{x\in\mathcal H:\ \langle x, x\rangle=0\}$.
Then the factor space $\mathcal H/\mathcal H_0$ with the inner product 
$\langle x+\mathcal H_0, y+\mathcal H_0\rangle=\langle x,y\rangle$ will 
be an inner product space.
Let $\mathcal H_+^{(a)}$ denote the resulting Hilbert space obtained 
by completion, and let $X_{T,+}$  
be the natural imbedding of $\mathcal H$ to $\mathcal H_+^{(a)}$,
$X_{T,+}x=x+\mathcal H_0$. Clearly, $X_{T,+}$ is a (linear, bounded) operator, that is, $X_{T, +} \in\mathcal L(\mathcal H,\mathcal H^{(a)}_+)$, 
and $\|X_{T,+}\|\leq \sup_{n\geq 0}\|T^n\|$. 
Furthermore, $\langle Tx,Ty\rangle=\langle x,y\rangle$ for every $x,y\in \mathcal H$.
Therefore, $T_1: x+\mathcal H_0\mapsto Tx+\mathcal H_0$ is a well-defined isometry on 
$\mathcal H/\mathcal H_0$. Denote by $T_+^{(a)}$ the continuous extension of $T_1$ to the space
$\mathcal H_+^{(a)}$. Clearly, $X_{T,+}$ realizes the relation $T\buildrel d \over \prec T_+^{(a)}$. If $T$ is not of class $C_{0\cdot}$, then 
 $T_+^{(a)}$ is not a zero operator. If the (nonzero) isometry $T_+^{(a)}$ is not a unitary operator, then $T$ has nontrivial  hyperinvariant subspaces
(see Introduction). 

The {\it unitary asymptote} 
$$\bigl(X_T \in\mathcal L(\mathcal H,\mathcal H^{(a)}),\ T^{(a)}\in\mathcal L(\mathcal H^{(a)})\bigr)$$ of $T$  is the minimal unitary extension of 
the  isometric asymptote of $T$. 

Suppose  that $T\in \mathcal L(\mathcal H)$ is a power bounded operator such that $T^{(a)}_+$ is a unitary operator, then, in fact, $T^{(a)}=T^{(a)}_+$ and
$X_T=X_{T, +}$. Furthermore, 
suppose that  $\mathcal M\in\operatorname{Lat}T$ is such that $T|_{\mathcal M}\approx S$. 
Then $X_T|_{\mathcal M}$ is bounded below, $X_T\mathcal M\in\operatorname{Lat}T^{(a)}$, and 
$T^{(a)}|_{X_T\mathcal M}\cong S$. Set 
$$\mathcal N=\bigvee_{n\geq 0}\bigl(T^{(a)}\bigr)^{-n}X_T\mathcal M.$$
Then $\mathcal N$ is a reducing subspace of $T^{(a)}$, and 
$T^{(a)}|_{\mathcal N}\cong U_{\mathbb T}$. Without loss of generality suppose that  $T^{(a)}|_{\mathcal N} = U_{\mathbb T}$, 
 $X_T\mathcal M=H^2$, $\mathcal N=L^2$, and 
$$\mathcal H^{(a)} = L^2\oplus\mathcal K=H^2\oplus(H^2_-\oplus\mathcal K),$$
where $\mathcal K$ is a reducing subspace of $T^{(a)}$.

Set $T_0=P_{\mathcal M^\perp}T|_{\mathcal M^\perp}$.
Operators $T$,  $T^{(a)}$ and $X_T$ have the forms
$$T\!=\!\begin{pmatrix} T|_{\mathcal M} & T_2 \cr \mathbb O & T_0\end{pmatrix}, \ 
T^{(a)}\! =\! \begin{pmatrix} S & (\cdot, \chi^{-1}\oplus 0)\chi^0 \cr \mathbb O & V\end{pmatrix},
 \text{ and } X_T\! =\! \begin{pmatrix} X_T|_{\mathcal M} & X_2 \cr \mathbb O & X_1 \end{pmatrix} $$
 with respect to the decompositions $\mathcal H = \mathcal M\oplus{\mathcal M}^\perp$ and ${\mathcal H}^{(a)} = H^2\oplus (H^2_-\oplus\mathcal K)$, where
$T_2$,  $X_2$, $X_1$, and $V$  are appropriate operators. Following \cite{fer}, set 
$$ X=\begin{pmatrix} {X_T}|_{\mathcal M} & X_2 \cr \mathbb O & I_{{\mathcal M}^\perp}\end{pmatrix}, \ \ \ X\in\mathcal L(\mathcal H, H^2\oplus{\mathcal M}^\perp).$$
Clearly, $X$ is invertible, and 
\begin{align}\label{E:5.1}T_1:=XTX^{-1} = \begin{pmatrix} S & \bigl(\cdot, X_1^\ast (\chi^{-1}\oplus 0)\bigr)\chi^0 \cr \mathbb O & T_0\end{pmatrix}. \end{align}
That is, $T\approx T_1$. Thus, one can consider the operator of the form \eqref{E:5.1} instead of $T$. Also, 
  $(I_{H^2}\oplus X_1) T_1 = T^{(a)}(I_{H^2}\oplus X_1)$, and 
  $\operatorname{clos} X_1 {\mathcal M}^\perp= H^2_-\oplus\mathcal K$. Set $X_0=P_{H^2_-\oplus\{0\}}X_1$. Clearly, 
$X_0\in\mathcal L({\mathcal M}^\perp, H^2_-)$, $\operatorname{clos} X_0 {\mathcal M}^\perp= H^2_-$, 
$X_1^\ast( \chi^{-1}\oplus 0) = X_0^\ast \chi^{-1}$, and $X_0T_0 = S_\ast X_0$.

\smallskip

The following proposition is a version of {\cite[Theorem 1]{fer}}.
\begin{proposition}\label{P:5.1}  Suppose  that $T_0\in\mathcal L(\mathcal H_0)$, $X_0\in\mathcal L(\mathcal H_0, H^2_-)$ and $X_0 T_0=S_\ast X_0$. 
Put $$T=\begin{pmatrix} S & (\cdot , X_0^\ast \chi^{-1})\chi^0 \cr \mathbb O & T_0\end{pmatrix}.$$
Then $(I_{H^2}\oplus X_0)T = U_{\mathbb T}(I_{H^2}\oplus X_0)$, and 
\begin{align}\label{E:5.2}  P_{H^2\oplus\{0\}}\varphi(T)x = P_+\varphi\cdot(X_0 x) \ \ \ \text{ for every } x\in\mathcal H_0 \end{align} 
and every polynomial $\varphi$.

 Furthermore, $T$ is power bounded if and only if $T_0$ is power bounded, 
$T$ is polynomially bounded if and only if $T_0$ is polynomially bounded, 
$T$ is a.c. polynomially bounded if and only if $T_0$ is a.c. polynomially bounded, 
and $T$ is similar to a contraction if and only if $T_0$ is similar to a contraction.

Moreover, if  $T_0$ is a power bounded operator of class $C_{0\cdot}$, and $\operatorname{clos}X_0\mathcal H_0=H^2_-$, then
$(I_{H^2}\oplus X_0, U_{\mathbb T})$ is the isometric asymptote of $T$, and $T$ is quasianalytic if and only if $\ker X_0=\{0\}$ and 
$$P_-L^2(\tau,m)\cap  X_0\mathcal H_0=\{0\}$$ for every Borel set $\tau\subset\mathbb T$ such that $m(\tau)<1$.\end{proposition}
\begin{proof} An easy computation shows that 
\begin{align}\label{E:5.3} P_{H^2\oplus\{0\}}T^nx = \sum_{k=0}^{n-1}(X_0 x, \chi^{k-n})\chi^k \ \ \ \text{ for every } x\in\mathcal H_0 \text{ and } n\geq 1. \end{align}
Equality \eqref{E:5.2} follows  easily from \eqref{E:5.3}. Furthermore,  equalities \eqref{E:5.3} imply that $T$ is power bounded if and only if 
$T_0$ is power bounded and there exists $C>0$ such that 
\begin{align}\label{E:5.4} \sum_{k=0}^{n-1}|(X_0 x,\chi^{k-n})|^2\leq C\|x\|^2 \ \ \ \text{ for every } x\in \mathcal H_0 \text{ and }  n\geq 1. \end{align}
Since
$$ \sum_{k=0}^{n-1}|(X_0 x,\chi^{k-n})|^2 = \sum_{k=1}^{n}|(X_0 x,\chi^{-k})|^2\leq \|X_0 x\|^2\leq \|X_0\|^2\|x\|^2,$$
 \eqref{E:5.4} is fulfilled.
Also, $\|P_+\varphi\cdot(X_0 x)\|\leq \|\varphi\|_\infty\|X_0\|\|x\|$ for every polynomial $\varphi$ and every $x\in \mathcal H_0$, 
and the conclusion on the polynomially boundedness is obtained. Moreover,  
the right part of \eqref{E:5.2} is defined for every $\varphi\in H^\infty$.  Therefore, 
if $T_0$ is a.c. polynomially bounded, then $\varphi(T)$ can be defined using \eqref{E:5.2}, and it is easy to check that 
the mapping 
$$\varphi\mapsto \varphi(T), \ \ H^\infty \to \mathcal L(H^2\oplus \mathcal H_0)$$ is a weak-$\ast$ continuous algebra-homomorphism. 
Now the conclusion on the a.c. polynomially boundedness follows from {\cite[Theorem 23]{ker5}}, see also \cite{mlak}.

The conclusion on the similarity to a contraction follows  from  {\cite[Corollary 4.2]{cass}}. 

Suppose that $T_0$ is a power bounded operator of class $C_{0\cdot}$. Equalities \eqref{E:5.3} imply that 
$$\lim_{n\to\infty}\big(T^n(h_1\oplus x_1),T^n(h_2\oplus x_2)\big)=(h_1,h_2)+(X_0x_1, X_0x_2)$$
for all $h_1$, $h_2\in H^2$, $x_1,x_2\in\mathcal H_0$. The conclusion on the isometric asymptote of $T$ follows from the latter equality and the results in \cite{ker1} (see the beginning of this section). The conclusion on the quasianalyticity is a straightforward consequence of the form of 
the  isometric asymptote of $T$. \end{proof}

The following theorem is a consequence of  the results in \cite{bercpr} (based on \cite{bour}) and \cite{mlak}, see also \cite{ker5}.

\begin{theorem}\label{T:5.2} Suppose that $R$ is a polynomially bounded operator. Then  $S\prec R$ if and only if $R$ is a.c., 
$R$ is not a $C_0$-operator (that is, there is no nonzero function $\varphi\in H^\infty$ such that $\varphi(R)=\mathbb O$), 
and $R$ is cyclic.\end{theorem}
\begin{proof} The ``only if" part is almost evident. We mention only that $R$ is a.c. by {\cite[Proposition 16]{ker5}}, see also \cite{mlak}.
To prove ``if" part, suppose that $R\in\mathcal L(\mathcal H)$ is a cyclic polynomially 
bounded operator with a cyclic vector $x_0\in \mathcal H$. By {\cite[Lemma 2.1]{bercpr}}, there exist a finite positive Borel measure $\mu$ 
on $\mathbb T$ and $X\in\mathcal L( P^2(\mu), \mathcal H)$ such that $XS_\mu=RX$ and $X\text{\bf 1}=x_0$. Here 
$P^2(\mu)$ is the closure of analytic polynomials in $L^2(\mu)$,  $S_\mu$ is 
the operator of multiplication by the independent variable in $P^2(\mu)$, and $\text{{\bf 1}}(\zeta)=\zeta$ for a.e. $\zeta\in\mathbb T$ 
(with respect to $\mu$). Since $x_0$ is cyclic for $R$, $\operatorname{clos}X P^2(\mu)= \mathcal H$. 

There exist a nonnegative function $w\in L^1(\mathbb T,m)$ and a finite positive Borel measure $\mu_s$ on $\mathbb T$, singular with respect to $m$, 
such that $\mu=wm+\mu_s$. By {\cite[Proposition III.12.3]{conw2}},  $P^2(\mu)=P^2(wm)\oplus L^2(\mu_s)$. Denote by $U_{\mu_s}$ 
the operator of multiplication by the independent variable in $L^2(\mu_s)$. Then  $U_{\mu_s}$ is a singular unitary operator. 
 We  have $X|_{L^2(\mu_s)}U_{\mu_s}=RX|_{L^2(\mu_s)}$. Now suppose that $R$ is a.c.. By 
{\cite[Proposition 15]{ker5}},  see also \cite{mlak}, $X|_{L^2(\mu_s)}=\mathbb O$. 
Set $X_1=X|_{P^2(wm)}$. We have  $\operatorname{clos}X_1 P^2(wm)= \mathcal H$ and $X_1S_{wm}=RX_1$. That is, $S_{wm}\buildrel d \over \prec R$.  
Recall that $S \prec S_{wm}$. Indeed, set  $\tau=\{\zeta\in\mathbb T\ :\ w(\zeta)=0\}$. 
If $m(\tau)>0$, then $ S_{wm}\cong U(\mathbb T\setminus\tau)$, and the relation $S \prec  U(\mathbb T\setminus\tau)$ is realized by the natural imbedding.
If $m(\tau)=0$ and $\log w\notin L^1(\mathbb T,m)$, then $ S_{wm}\cong U_\mathbb T$, and the relation $S \prec  U_\mathbb T$ is realized by the 
operator of multiplication by  some function $\psi\in L^\infty(\mathbb T,m)$ such that $\psi\neq 0$ a.e.  and 
$\log |\psi|\notin L^1(\mathbb T,m)$. 
If $\log w\in L^1(\mathbb T,m)$, then $ S_{wm}\cong S$ (see, for example, {\cite[Ch. III.12, VII.10]{conw2}}).

Thus, if $R\in\mathcal L(\mathcal H)$ is a cyclic a.c. polynomially bounded operator, then $S\buildrel d \over \prec R$.  
Denote by $Y\in\mathcal L( H^2, \mathcal H)$  an operator which realizes this relation. If $\ker Y\neq \{0\}$, then there exists 
an inner function $\theta$ such that $\ker Y=\theta H^2$ (because $\ker Y\in\operatorname{Lat}S$). Set $K_\theta=H^2\ominus\theta H^2$ and  
$Y_0=Y|_{K_\theta}$. Then $Y_0$ realizes the relation $P_{K_\theta}S|_{K_\theta}\prec R$. Since $\theta(P_{K_\theta}S|_{K_\theta})=\mathbb O$, 
we conclude that $\theta(R)=\mathbb O$. Thus, if  $R$ is a cyclic a.c. polynomially bounded operator, and $R$ is not a $C_0$-operator, 
then $S\prec R$.  
\end{proof}

\begin{corollary}\label{C:5.3} Suppose that $T_0$ is a polynomially bounded operator. Then  
there exists  a quasiaffinity $X_0$ such that $X_0T_0=S_\ast X_0$ if and only if $T_0$ is a.c., $T_0$ is not a $C_0$-operator, 
and $T_0^\ast$ is cyclic.\end{corollary}
\begin{proof} Clearly,  $(S_\ast)^\ast\cong S$. It remains to apply Theorem~\ref{T:5.2} to $T_0^\ast$. \end{proof}

\begin{proposition}\label{P:5.4} Suppose  that $T_0$, $X_0$ and $T$ are as in Proposition~\ref{P:5.1}, 
 $T_0$ is a.c. polynomially bounded, $\operatorname{clos}X_0\mathcal H_0=H^2_-$,  and $\theta\in H^\infty$ is an inner function.
Then $\ker\theta(T)^\ast\neq\{0\}$ if and only if there exists $x\in\mathcal H_0$ such that $x\notin X_0^\ast H^2_-$ and
$\theta(T_0)^\ast x \in  X_0^\ast H^2_-$.\end{proposition}
\begin{proof} Set $K_\theta = H^2\ominus\theta H^2$, then 
$K_\theta = \theta\overline\chi\overline{K_\theta} = H^2\cap\theta H^2_-$.
 Note that $\widetilde\theta((S_\ast)^\ast)f = \overline\theta f$ for every $f\in H^2_-$.  
Suppose that  $x\in\mathcal H_0$ is such that $\theta(T_0)^\ast x \in  X_0^\ast H^2_-$. Then there exist $f\in H^2_-$ and 
$g\in K_\theta$ such that 
\begin{align*}\theta(T_0)^\ast x & = X_0^\ast(\overline\chi\overline g + \overline\theta f) = 
X_0^\ast\overline\chi\overline g +  X_0^\ast \widetilde\theta((S_\ast)^\ast)f \\
& =  X_0^\ast\overline\chi\overline g +  \widetilde\theta(T_0^\ast)X_0^\ast f = X_0^\ast\overline\chi\overline g + \theta(T_0)^\ast X_0^\ast f. \end{align*}
From the definition of $T$ and \eqref{E:5.2} applied to $\theta$ we have that 
$$\theta(T)^\ast\bigl(\theta\overline\chi\overline g\oplus(X_0^\ast f-x)\bigr)= P_+\overline\theta\theta\overline\chi\overline g\oplus
\bigl(X_0^\ast P_-\overline\theta\theta\overline\chi\overline g + \theta(T_0)^\ast(X_0^\ast f-x)\bigr) =0.$$
If   $x\notin X_0^\ast H^2_-$, then $X_0^\ast f-x\neq 0$. Therefore, $\ker\theta(T)^\ast\neq\{0\}$.

Now let $h\in H^2$ and $x\in\mathcal H_0$ be such that $h\oplus x\neq 0$ and $\theta(T)^\ast(h\oplus x)=0$. 
Then $\theta(S)^\ast h=P_+\overline\theta h =0$, 
therefore, $h\in K_\theta$, and $h=\theta\overline\chi\overline g$ for some $g\in K_\theta$. 
Furthermore, $\theta(T_0)^\ast x =-X_0^\ast \overline\chi\overline g$, and we conclude that $\theta(T_0)^\ast x\in X_0^\ast  H^2_-$.
If there exists $f\in H^2_-$ such that $x = X_0^\ast f$, then 
$$\theta(T_0)^\ast x = \theta(T_0)^\ast X_0^\ast f =X_0^\ast\widetilde\theta((S_\ast)^\ast)f = X_0^\ast\overline\theta f.$$
Since $\ker X_0^\ast  = \{0\}$, we conclude that $\overline\theta f = -\overline\chi\overline g$. Taking into account that 
 $g\in K_\theta$ and $f\in H^2_-$, we conclude that $g=0$ and $f=0$, consequently, $h=0$ and $x=0$, a contradiction. Therefore, 
$x\notin X_0^\ast H^2_-$.   \end{proof}

\begin{proposition}\label{P:5.5} Suppose   that $T_0\in\mathcal L(\mathcal H_0)$, $X_0\in\mathcal L(\mathcal H_0, H^2_-)$, $X_0 T_0=S_\ast X_0$, 
\begin{align}\label{E:5.5}\bigvee_{n\geq 0}(T_0^\ast)^nX_0^\ast\chi^{-1} = \mathcal H_0, \end{align} and  there is a sequence $\{y_n\}_{n=0}^\infty\subset\mathcal H_0$
 such that \begin{equation}\label{E:5.6}\begin{gathered}\bigvee_{n\geq 0}y_n = \mathcal H_0 \ \ \text{ and } \\ \bigl((T_0^\ast)^nX_0^\ast\chi^{-1},y_k\bigr)=0, \ \ n\neq k, \ \  
\bigl((T_0^\ast)^nX_0^\ast\chi^{-1},y_n\bigr)=1, \ \ \ n,k\geq0.\end{gathered} \end{equation}
Furthermore, suppose  that $T_0$ is  a.c. polynomially bounded,  $\theta\in H^\infty$ is a singular inner function, 
 and $x\in\mathcal H_0$. 
Then $x\notin X_0^\ast H^2_-$ and $\theta(T_0)^\ast x \in X_0^\ast H^2_-$ if and only if there exists $f\in H^2$ such that 
$f\notin\theta H^2$ and \begin{align}\label{E:5.7}(x, y_n)=\overline{\widehat{\frac{f}{\theta}(}n)} \ \ \ \text{ for all } n\geq 0.\end{align}\end{proposition}
\begin{proof} It is easy to see that $X_0y_n=\chi^{-n-1}$ for all $n\geq 0$, therefore,
\begin{align}\label{E:5.8}   (X_0^\ast\chi^{-1}\overline f, y_n)= \overline{{\widehat f}(n)}\end{align}  
for every $f\in H^2$ and $n\geq 0$. Also, 
\begin{align}\label{E:5.9}  (\varphi(T_0)^\ast (T_0^\ast)^nX_0^\ast\chi^{-1},y_k) = \overline{{\widehat \varphi}(k-n)} \end{align}  
for every $\varphi\in H^\infty$, $n$, $k\geq 0$. It  follows from   \eqref{E:5.5} that for every $x\in\mathcal H_0$ 
there exist a sequence $\{N_j\}_j$ of positive integers
 and a sequence of  families  $\{a_{jn}\}_{n=0}^{N_j}$ of complex numbers such that 
\begin{align}\label{E:5.10}  x=\lim_j\sum_{n=0}^{N_j}a_{jn}(T_0^\ast)^nX_0^\ast\chi^{-1}. \end{align}  
From \eqref{E:5.6} and \eqref{E:5.10} we obtain that\begin{align}\label{E:5.11}  (x,y_n)=\lim_ja_{jn}.\end{align}  
Let $\varphi\in H^\infty$. We infer from \eqref{E:5.9}, \eqref{E:5.10}, and \eqref{E:5.11} that 
\begin{equation}\label{E:5.12}\begin{aligned}(\varphi(T_0)^\ast x,y_k)& = \lim_j\sum_{n=0}^{N_j}a_{jn}\bigl(\varphi(T_0)^\ast(T_0^\ast)^nX_0^\ast\chi^{-1},y_k\bigr) 
=\lim_j\sum_{n=0}^{N_j}a_{jn}\overline{{\widehat \varphi}(k-n)}\\
&=\lim_j\sum_{n=0}^k a_{jn}\overline{{\widehat \varphi}(k-n)} = \sum_{n=0}^k (x,y_n)\overline{{\widehat \varphi}(k-n)}. \end{aligned}\end{equation}

Now suppose that $\theta(T_0)^\ast x \in  X_0^\ast H^2_-$, that is, there exists $f\in H^2$ such that 
$\theta(T_0)^\ast x = X_0^\ast \chi^{-1}\overline f$. From \eqref{E:5.8} and \eqref{E:5.12} we conclude that 
\begin{align}\label{E:5.13}  (\theta(T_0)^\ast x, y_k) = \sum_{n=0}^k (x,y_n)\overline{{\widehat \theta}(k-n)} = \overline{{\widehat f}(k)} 
\ \ \text{ for all } k\geq 0. \end{align}
Equalities \eqref{E:5.7} follow from equalities \eqref{E:5.13}, because $\widehat\theta(0)\neq 0$. If $f=\theta h$ for some $h\in H^2$, then 
it  follows from \eqref{E:5.7} and \eqref{E:5.8} that $x=X_0^\ast\chi^{-1}\overline h$. Part ``only if" is proved.

Now suppose that \eqref{E:5.7} is fulfilled for some $f\in H^2$. From the equality 
\begin{align}\label{E:5.14}\overline{{\widehat f}(k)} = \sum_{n=0}^k (x,y_n)\overline{{\widehat \theta}(k-n)}, 
\end{align} \eqref{E:5.8}, and \eqref{E:5.12} we conclude that 
$\theta(T_0)^\ast x =  X_0^\ast \chi^{-1}\overline f$. If  $x=X_0^\ast\chi^{-1}\overline h$ for some $h\in H^2$, 
 then it follows from \eqref{E:5.8} that  $(x,y_n)=\overline{{\widehat h}(n)}$ for all $n\geq 0$. These equalities and 
\eqref{E:5.14} imply that $\widetilde f(z) =\widetilde\theta(z)\widetilde h(z)$  for all $z\in\mathbb D$, therefore, $f=\theta h$, that is, $f\in\theta H^2$.
Part ``if" is proved.  \end{proof}

Applying Propositions~\ref{P:5.4} and \ref{P:5.5} to a dissymmetric weighted shift (see Sec. 2 for definitions and references), we obtain the following known result.

 \begin{corollary} \label{C:5.6} Suppose  that $\omega$ is a dissymmetric weight, and  $\theta\in H^\infty$ is a singular inner function.
Then $\operatorname{ran}\theta(S_\omega)$ is not dense if and only if there exists $f\in H^2$ such that $f\not\in\theta H^2$ and 
$$\sum_{n=0}^\infty \frac{1}{\omega(-1-n)^2}\Bigl|\widehat{\frac{f}{\theta}(}n)\Bigr|^2<\infty. $$\end{corollary}

\begin{remark} \label{R:5.7} If $T_0$ and $X_0$ satisfy  conditions \eqref{E:5.5} and \eqref{E:5.6}, $\ker X_0=\{0\}$, and $T$ is constructed by $T_0$ and $X_0$ as in 
Proposition~\ref{P:5.1}, then it is easy to see that $Ty_{n+1}=y_n$, $n\geq 0$. Let $T_0$ be power bounded. By Proposition~\ref{P:5.1}, $T$ is power bounded 
too. If 
$$\sum_{n=0}^\infty\frac{\log\|y_n\|}{n^2+1}<\infty,$$
then $T$ is not quasianalytic, and, consequently, has nontrivial hyperinvariant subspaces (see, for example, 
{\cite[Theorem XII.8.1]{beau}},  \cite{kell2}, \cite{ker2}).\end{remark}

 \begin{theorem}  \label{T:5.8} Suppose that  $T_0$, $X_0$ and $T$ are as in Proposition~\ref{P:5.1}, and there exists a sequence of positive numbers 
$\{\beta_n\}_{n=0}^\infty$ such that $\beta_n\to\infty$ and 
$$\sum_{n=0}^\infty\bigl|\bigl(x,(T_0^\ast)^n X_0^\ast\chi^{-1}\bigr)\bigr|^2\beta_n^2<\infty \ \ \ \text{ for every } x\in \mathcal H_0. $$
  Then there exist a dissymmetric weight $\omega$ and $Y_0\in\mathcal L(\mathcal H_0, \ell^2_{\omega -})$ such that  
$(J\oplus Y_0)T = S_\omega(J\oplus Y_0)$, where $J\in\mathcal L(H^2,\ell^2_+)$ is a unitary operator acting by the formula 
$Jh=\{\widehat{h}(n)\}_{n\geq 0}$, $h\in H^2$. 

Set $\mathcal H=H^2\oplus\mathcal H_0$ and $Y=J\oplus Y_0$.  If $\operatorname{clos}T\mathcal H=\mathcal H$, 
then $\operatorname{clos}Y\mathcal H=\ell^2_\omega$, that is, $T\buildrel d \over \prec S_\omega$.\end{theorem}
\begin{proof} By Lemma~\ref{L:2.2}(ii), there exists a dissymmetric weight $\omega$ such that
$\omega(-n-1)\leq\beta_n$ for sufficiently large $n$. Define $Y_0$ by the formulas 
\begin{align*}(Y_0x)(n)=0,  \ n\geq 0, \ \ (Y_0x)(n)=\bigl(x, (T_0^\ast)^{-n-1}X_0^\ast\chi^{-1}\bigr), \ \ n\leq -1, \ \ x\in\mathcal H_0. \end{align*} Then $Y_0$ is bounded by the closed graph theorem.
The  intertwining relation $YT=S_\omega Y$ follows from an easy calculation.
 Since $S_\omega$ is invertible, we have $Y=S_\omega^{-1}YT$. Now suppose that $\operatorname{clos}T\mathcal H=\mathcal H$. 
We have $$S_\omega^{-1}\operatorname{clos}Y\mathcal H=S_\omega^{-1}\operatorname{clos}Y\operatorname{clos}T\mathcal H
\subset\operatorname{clos}S_\omega^{-1}YT\mathcal H=\operatorname{clos}Y\mathcal H.$$ It means that 
$\operatorname{clos}Y\mathcal H\in\operatorname{Lat}S_\omega^{-1}$. Furthermore, $\ell^2_+\subset\operatorname{clos}Y\mathcal H$, and we conclude that 
$\operatorname{clos}Y\mathcal H=\ell^2_\omega$.   \end{proof}

\begin{corollary} \label{C:5.9} Suppose that in the conditions of Theorem~\ref{T:5.8}, $T$ is a.c. polynomially bounded and 
$\operatorname{clos}T\mathcal H=\mathcal H$. Then there exists a singular inner function $\theta$ 
such that $m(\operatorname{supp}\mu_\theta)=0$ and $\operatorname{clos}\theta(T)\mathcal H\neq\mathcal H$.\end{corollary}
\begin{proof} This is a straightforward consequence of Theorem~\ref{T:5.8}, Lemma~\ref{L:1.1} and 
Theorem~\ref{T:2.1}. \end{proof}

Recall that $A^+(\mathbb T)$ is the  algebra of analytic functions with absolutely summable Taylor coefficients 
(see Sec. 1.2).
If $T$ is power bounded, that for every $\varphi\in A^+(\mathbb T)$ the operator $\varphi(T)$ is defined by the formula
$$\varphi(T)=\sum_{n=0}^\infty\widehat\varphi(n)T^n.$$ Let $\theta$ be a singular inner function. Then the condition 
$m(\operatorname{supp}\mu_\theta)=0$ is necessary for the existence 
$\psi\in H^\infty$ such that $\theta\psi\in A^+(\mathbb T)$, but is not sufficient (\cite{carl}, \cite{kaka}).
Namely, there exist singular inner functions $\theta$ such that $m(\operatorname{supp}\mu_\theta)=0$,
but there is no function $\varphi\in A^+(\mathbb T)$ such that $\varphi =0$ on  
$\operatorname{supp}\mu_\theta$ and $\varphi\not\equiv 0$. For such $\theta$, there is no (nonzero) $\psi\in H^\infty$ 
such that $\theta\psi\in A^+(\mathbb T)$. On the other hand, if $m(\operatorname{supp}\mu_\theta)=0$ and 
$\operatorname{supp}\mu_\theta$ is a  Carleson set (see Definition~\ref{D:2.3}), 
then there exists an outer function $\psi\in A^+(\mathbb T)$  such that $\theta\psi\in A^+(\mathbb T)$ \cite{tawi}. 

\begin{corollary} \label{C:5.10} Suppose that $T\in\mathcal L(\mathcal H)$,   $\{\beta_n\}_{n=0}^\infty$ is a sequence
 of positive numbers, $T$ and  $\{\beta_n\}_{n=0}^\infty$ satisfy the conditions of Theorem~\ref{T:5.8},
 $T$ is power bounded,  and $\operatorname{clos}T\mathcal H=\mathcal H$. 
Furthermore, suppose that there exists  a sequence of positive numbers $\{w_n\}_{n\geq 1}$ satisfying 
the conditions of Theorem~\ref{T:2.4} and  
such that $w_{n+1}\leq\beta_n$ for sufficiently large $n$. 
Then there exist a singular inner function $\theta$ and an outer function $\psi\in A^+(\mathbb T)$ such that 
$\theta\psi\in A^+(\mathbb T)$ and  $\operatorname{clos}(\theta\psi)(T)\mathcal H\neq\mathcal H$.\end{corollary}
\begin{proof}
Let $0<s<1$. By Theorem~\ref{T:2.4} and Lemma~\ref{L:2.5}(ii), there exist a singular inner function $\theta$ and a constant $C>0$ such that $\operatorname{supp}\mu_\theta$ is a Carleson set,  
$m(\operatorname{supp}\mu_\theta)=0$, and 
$$ \Bigl|\widehat{\frac{1}{\theta}}(n-1)\Bigr|\leq C w_n^s \ \ \ \text{ for all } n\geq 1.$$
Let $\omega$ be the  dissymmetric weight constructed by Lemma~\ref{L:2.5}(i) applied to  $\{w_n\}_{n\geq 1}$. We have
$$\sum_{n\geq 0}\frac{1}{\omega(-1-n)^2} \Bigl|\widehat{\frac{1}{\theta}}(n)\Bigr|^2\leq
\sum_{n\geq 0}\frac{C^2 w_{n+1}^{2s}}{w_{n+1}^2}=C^2\sum_{n\geq 0}\frac{1}{w_{n+1}^{2(1-s)}}<\infty$$
by  Lemma~\ref{L:2.5}(ii). Thus, $\theta$ and $\omega$ satisfy \eqref{E:2.3}. By Theorem~\ref{T:2.1}, $\operatorname{ran} \theta(S_\omega)$ is not dense. Consequently, 
$\operatorname{ran} (\theta\psi)(S_\omega)$ is not dense for any $\psi\in H^\infty$. 

We prove that $T\buildrel d \over \prec S_\omega$ exactly as in the proof of Theorem~\ref{T:5.8}. By \cite{tawi}, 
there exists an outer function $\psi\in A^+(\mathbb T)$  such that $\theta\psi\in A^+(\mathbb T)$. Set $\varphi= \theta\psi$ and apply the evident analog of Lemma~\ref{L:1.1}  
to the power bounded operators $T$ and $S_\omega$ and $\varphi\in A^+(\mathbb T)$.   \end{proof}
 
\begin{remark}\label{R:5.12}  If $T\in\mathcal L(\mathcal H)$ is an a.c. contraction and $T \buildrel d \over \prec U_{\mathbb T}$, then 
$$\bigvee\{\mathcal M \ :\ \mathcal M\in\operatorname{Lat}T, \ T|_{\mathcal M}\approx S\}=\mathcal H$$
(\cite{ker3} or {\cite[Theorem IX.3.6]{sznagy}}). Indeed, let $(X,U)$ be the unitary asymptote of $T$, then 
$U$ is an a.c. unitary operator. Denote by $Y$ the operator which realizes the relation $T \buildrel d \over \prec U_{\mathbb T}$. 
Then there exists 
an operator $Z$ such that $Y=ZX$ and $ZU=U_{\mathbb T}Z$. Since $\operatorname{ran}Y$ is dense in $L^2$, $\operatorname{ran}Z$ is dense in $L^2$, 
too. Therefore, $U$ has a reducing subspace such that the restriction of $U$ on this subspace is unitarily equivalent to  $U_{\mathbb T}$. Thus, 
\cite{ker3} or {\cite[Theorem IX.3.6]{sznagy}} can be applied to $T$, and the needed conclusion is obtained. \end{remark}
 
\begin{remark}\label{R:5.13} It is proved in \cite{kert} that if every contraction $T$ whose unitary asymptote is $U_\mathbb T$ has a nontrivial hyperinvariant 
subspace, then every contraction which has a nonzero {\it cyclic} unitary asymptote (and is not a multiple of the identity operator) has 
a nontrivial hyperinvariant subspace.\end{remark}

\section{Sufficient conditions for the  existence of hyperinvariant subspaces for operators intertwined with a.c. unitaries}

Let $\varphi\in H^\infty$, and let $\mu$ be a positive finite Borel measure on $\mathbb T$. For $\xi\in\mathbb T$ set
$$\varphi_\xi(z)=\varphi(\xi z), \ \ z\in\mathbb D, \ \ \text{ and } \ \ \mu_\xi(\tau)=\mu(\xi\tau), \ \ \tau\subset \mathbb T, $$ where $\xi\tau=\{\xi\zeta:\ \zeta\in\tau\}$.
Let $\theta$ be a singular inner function. Since $\theta$ has no zeros in $\mathbb D$, the function $1/\theta$ is analytic in $\mathbb D$, indeed,   
$$\frac{1}{\theta(z)}=\sum_{n=0}^\infty \widehat{\frac{1}{\theta}}(n) z^n, \ \ z\in\mathbb D.$$
From the equality $\frac{1}{\theta(z)}\theta(z) = 1$, $z\in\mathbb D$, we have 
\begin{align}\label{E:6.1} \sum_{k=0}^n \widehat{\frac{1}{\theta}}(k)\widehat\theta(n-k)=0, \ \ n\geq 1, \ \ \ 
\widehat{\frac{1}{\theta}}(0)\widehat\theta(0) = 1. \end{align}
Recall that $\mu_\theta$ and $\widetilde \theta$ are defined in Sec. 1.2. We have 
$$(\theta_\xi)\widetilde{\ } =(\widetilde \theta)_{\overline\xi}, \ \ \ \mu_{\theta_\xi}=(\mu_\theta)_\xi, \ \ \
\widehat{\frac{1}{\theta_\xi}}(n)=\widehat{\frac{1}{\theta}}(n)\xi^n, \ \ n\geq 0.$$

\begin{lemma}\label{L:6.1} Suppose $\mu$ is a positive finite Borel singular measure
on $\mathbb T$, and $\nu$ is a positive Borel measure on $\mathbb T$. Let
$\Xi=\{\xi\in\mathbb T :\ \nu(\tau)\geq\mu_\xi(\tau)$ for every Borel set $\tau\subset \mathbb T\}$. If there exists nonempty open  arc
 $\mathcal J\subset\operatorname{clos}\Xi$, then $\nu(\mathbb T)=\infty$.\end{lemma}
\begin{proof} For $\zeta\in\mathbb T$ and $t>0$ set $\mathcal J_{\zeta,t}=\{\zeta e^{is} :\ |s|\leq t\}$. Since $\mu$ is singular, there exists  $\zeta_0\in\mathbb T$ such that 
\begin{align*}\lim_{t\to 0}\frac{1}{t}\mu(\mathcal J_{\zeta_0,t})=\infty \end{align*}  
({\cite[(II.6.3)]{garn}} or {\cite[Proposition III.4.1]{conw2}}).

For $n\in\mathbb N$ set $t_n=\pi m(\mathcal J)/2n$. There exists a family $\{\xi_k\}_{k=1}^n\subset\Xi$ such that 
$\overline\xi_k\mathcal J_{\zeta_0,t_n}\cap\overline\xi_l\mathcal J_{\zeta_0,t_n}=\emptyset$ for $k\neq l$. We have 
\begin{align*}\nu(\mathbb T)\geq \sum_{k=1}^n \nu(\overline\xi_k\mathcal J_{\zeta_0,t_n})&\geq 
\sum_{k=1}^n \mu(\xi_k\overline\xi_k\mathcal J_{\zeta_0,t_n}) = 
n\mu(\mathcal J_{\zeta_0,t_n})\\
& = n t_n \frac{1}{t_n}\mu(\mathcal J_{\zeta_0,t_n}) = 
\frac{\pi m(\mathcal J)}{2}\frac{1}{t_n}\mu(\mathcal J_{\zeta_0,t_n})\to\infty\end{align*}
when $n\to\infty$.  \end{proof}

\begin{corollary}\label{C:6.2} Suppose $\vartheta$ and $\theta$ are singular inner functions. Then the closure of the set of $\xi\in\mathbb T$ 
such that $\overline{\theta_\xi}\vartheta\in H^\infty$
does not contain nonempty open arcs. \end{corollary}
\begin{proof} Set $\nu=\mu_\vartheta$ and $\mu=\mu_\theta$. If $\xi\in\mathbb T$ is such that $\overline{\theta_\xi}\vartheta\in H^\infty$, 
then $\nu(\tau)\geq\mu_\xi(\tau)$ for every Borel set $\tau\subset \mathbb T$. Since $\nu(\mathbb T)<\infty$, the conclusion of the corollary 
follows from Lemma~\ref{L:6.1}. \end{proof}

\begin{lemma}\label{L:6.3} Let $\varphi\in H^\infty$, and let $f\in L^1(\mathbb T,m)$. Put 
$$(\varphi\ast f)(\xi)=\int_{\mathbb T}\varphi(\xi\zeta)f(\zeta)\text{{\rm d}}m(\zeta), \ \  \xi\in \mathbb D.$$
Then $\varphi\ast f$ is analytic in $\mathbb D$, continuous in $\operatorname{clos}\mathbb D$, 
$(\varphi\ast f)(\xi)=\int_{\mathbb T}\varphi_\xi f\text{{\rm d}}m$, $\xi\in\mathbb T$, and
\begin{align}\label{E:6.3}\widehat{\varphi\ast f}(n)=\widehat\varphi(n)\widehat f(-n),  \ \ \ n\geq 0.\end{align}
\end{lemma}
\begin{proof} It is easy to see that  $\varphi\ast f\in H^\infty$, and $\varphi\ast f$ has nontangential boundary value at every $\xi\in\mathbb T$, which is equal to 
 $\int_{\mathbb T}\varphi_\xi f\text{{\rm d}}m$.

The function $\xi\mapsto\int_{\mathbb T}\varphi_\xi f\text{{\rm d}}m$ is continuous on $\mathbb T$. Indeed, let $\xi$,
$\xi_0\in\mathbb T$. Then
$$\int_{\mathbb T}\varphi_{\xi_0} f\text{{\rm d}}m -\int_{\mathbb T}\varphi_\xi f\text{{\rm d}}m =
\int_{\mathbb T}\varphi(\zeta)(f(\zeta\overline{\xi_0})-f(\zeta\overline\xi))\text{{\rm d}}m(\zeta).$$
Therefore, 
\begin{equation*}\begin{aligned}\Bigl|\int_{\mathbb T}\varphi_{\xi_0} f\text{{\rm d}}m -\int_{\mathbb T}\varphi_\xi f\text{{\rm d}}m\Bigr|\leq
\|\varphi\|_\infty\int_{\mathbb T}|f(\zeta\overline\xi_0)-&f(\zeta\overline\xi)|\text{{\rm d}}m(\zeta)\to 0 \\ &\text{ when } \xi\to\xi_0,
\end{aligned}\end{equation*}
see, for example, {\cite[(I.2.1)]{katz}}. 

 Every function from $H^\infty$ is the convolution of its boundary values with the Poisson kernel. The  convolution of a continuous function on $\mathbb T$ with the Poisson kernel  is continuous in $\operatorname{clos}\mathbb D$. The conclusion of the lemma follows from these facts. Formula \eqref{E:6.3} for the Fourier coefficients of convolution can be easy 
checked straightforward in the conditions of the lemma.\end{proof}

\begin{lemma}\label{L:6.4} Suppose that  $T\in\mathcal L(\mathcal H)$, $\tau\subset \mathbb T$ is a Borel set, $m(\tau)>0$, 
$X\in\mathcal L(\mathcal H, L^2(\tau,m))$, $XT=U(\tau)X$, and $\operatorname{clos}X\mathcal H=L^2(\tau,m)$.
Suppose that  $\theta$ is a singular inner function, $0\not\equiv g\in L^2(\tau,m)$,  and
 $$\sum_{n=0}^\infty\Bigl |\widehat{\frac{1}{\theta}}(n)\Bigr|\|T^{\ast n}X^\ast g \|<\infty.$$ 
For every $\xi\in\mathbb T$ put 
$$u_\xi=\sum_{n=0}^\infty \widehat{\frac{1}{\theta_\xi}}(n)T^{\ast n}X^\ast g \ \ \ 
\text{ and } \ \ \ v_\xi=X^\ast (\theta_\xi)\widetilde{\ }g.$$
Let $\Xi=\{\xi\in\mathbb T :\ u_\xi= v_\xi\}$. Then $\Xi$ is closed and $\Xi\neq\mathbb T$.\end{lemma}
\begin{proof}
For $x\in\mathcal H$ and  integers $n$, $k\geq 0$ put $\psi_x=Xx$, $f_{x,k}=\overline\chi^k\overline{\psi_x} g$,  and
\begin{align}\label{E:6.4}a_{x,k,n}=\widehat{\frac{1}{\theta}}(n)(T^{\ast n}X^\ast g,T^kx) = \widehat{\frac{1}{\theta}}(n)\widehat{\overline{\psi_x} g}(n+k).\end{align}
Then $\sum_{n\geq 0}|a_{x,k,n}|<\infty$, \begin{align}\label{E:6.5} (u_\xi, T^kx)=\sum_{n\geq 0}a_{x,k,n}\xi^n, 
\end{align}
 and
\begin{align}\label{E:6.6}(v_\xi, T^kx)=\int_{\mathbb T}(\theta_\xi)\widetilde{\ }f_{x,k}\text{{\rm d}}m=(\widetilde\theta\ast f_{x,k})(\overline\xi). \end{align}
for every $\xi\in\mathbb T$.   

 By \eqref{E:6.5} and \eqref{E:6.6}, $\xi\in\Xi$ if and only if 
 \begin{align}\label{E:6.7}\sum_{n\geq 0}a_{x,k,n}\xi^n = (\widetilde\theta\ast f_{x,k})(\overline\xi)\ \ \ \text{ for all } x\in\mathcal H\ \text{ and } k\geq 0.\end{align}  
Since the functions on both sides of \eqref{E:6.7} are continuous on $\mathbb T$ (see the estimate before \eqref{E:6.5} and Lemma~\ref{L:6.3}), $\Xi$ is closed.
 
Let us assume that $\Xi=\mathbb T$. Then \eqref{E:6.7} implies that 
\begin{equation}\label{E:6.8}\begin{aligned} a_{x,k,n}=0, \ \ \ \widehat{\widetilde\theta\ast f_{x,k}}(n)=0, \ \ \text{ and } \ \ & a_{x,k,0}=\widehat{\widetilde\theta\ast f_{x,k}}(0) \\
&\text{ for every } \ \  n\geq 1, \ \ k\geq 0.\end{aligned}\end{equation}
By \eqref{E:6.3}, \eqref{E:6.4}, and \eqref{E:6.8}, 
$\overline{\widehat\theta(0)}\widehat{\overline{\psi_x} g}(k)=\widehat{\frac{1}{\theta}}(0)\widehat{\overline{\psi_x} g}(k)$. 
If $\widehat{\overline{\psi_x} g}(k)\neq 0$ for some $k\geq 0$, then $\overline{\widehat\theta(0)}=\widehat{\frac{1}{\theta}}(0)$, 
which contradicts to \eqref{E:6.1}, because $|\widehat\theta(0)|<1$. Thus,   
\begin{align}\label{E:6.9}\widehat{\overline{\psi_x} g}(k)=0 \ \ \text{ for all } \ \ k\geq 0. \end{align}
From \eqref{E:6.8} and  \eqref{E:6.9} we conclude that $(u_\xi, T^kx)=0$ for all $k\geq 0$ and $\xi\in\mathbb T$. 
Therefore, $(v_\xi, T^kx)=0$ for all $k\geq 0$ and $\xi\in\mathbb T$.  Taking into account \eqref{E:6.6}, we conclude that 
 $\overline{(\theta_\xi)\widetilde{\ }}\overline\chi\overline g \psi_x\in H^1$ for every  $\xi\in\mathbb T$. Thus, 
$\overline\chi\overline g \psi_x\in (\theta_\xi)\widetilde{\ }H^1$ for every  $\xi\in\mathbb T$. 
By Corollary~\ref{C:6.2}, 
$\overline g \psi_x\equiv 0$. In particular, $(Xx,g)=0$.

Since $x\in\mathcal H$ is arbitrary and $\operatorname{clos}X\mathcal H=L^2(\tau,m)$, we obtain that $g\equiv 0$, a contradiction. 
Thus, $\Xi\neq\mathbb T$.    \end{proof}

\begin{lemma}\label{L:6.5} Suppose that $S\in\mathcal L(H^2)$ is the unilateral shift of
 multiplicity $1$,
$H^2_\infty=\oplus_{n=1}^\infty H^2$, $S_\infty=\oplus_{n=1}^\infty S$. Suppose that 
$\theta\in H^\infty$ is a singular inner function, $h\in H^2_\infty$ is such that
$$\sum_{n=0}^\infty \Bigl|\widehat{\frac{1}{\theta}}(n)\Bigr|\|S_\infty^{\ast n} h\|<\infty, \  \ \text{ and } \ \ 
u = \sum_{n=0}^\infty \widehat{\frac{1}{\theta}}(n)S_\infty^{\ast n} h.$$ Then $\theta(S_\infty^\ast)u=h$.\end{lemma}
\begin{proof} We have $h=\{h_k\}_{k=1}^\infty$, where $h_k\in H^2$.
Since $S_\infty^{\ast n} h=\{S^{\ast n}h_k\}_{k=1}^\infty$, we have 
$\sum_{n=0}^\infty \bigl|\widehat{\frac{1}{\theta}}(n)\bigr|\|S^{\ast n} h_k\|<\infty$ for every $k\geq 1$. Set 
$$u_k= \sum_{n=0}^\infty \widehat{\frac{1}{\theta}}(n)S^{\ast n} h_k,$$ then $u=\{u_k\}_{k=1}^\infty$,
and $\theta(S_\infty^\ast)u= \{\theta(S^\ast) u_k\}_{k=1}^\infty$. 
Thus, it is sufficient to prove that $\theta(S^\ast) u_k = h_k$ for every $k\geq 1$.

For fixed $k$ put $f=h_k$ and $g=u_k$. We have for $j\geq 0$ 
\begin{align*}(\theta(S^\ast)g, \chi^j)=(g, \widetilde\theta(S)\chi^j)  & = 
\Bigl(\sum_{n=0}^\infty \widehat{\frac{1}{\theta}}(n)S^{\ast n}f, \widetilde\theta(S)\chi^j\Bigr) \\ 
& = 
 \sum_{n=0}^\infty \widehat{\frac{1}{\theta}}(n)(S^{\ast n}f, \widetilde\theta(S)\chi^j),\end{align*} 
and
$(S^{\ast n}f, \widetilde\theta(S)\chi^j)=(f,\widetilde\theta\chi^{j+n})
=\sum_{l=n+j}^\infty \widehat f(l)\widehat\theta(l-n-j)$.
Thus, $$(\theta(S^\ast)g, \chi^j) = \sum_{n=0}^\infty \widehat{\frac{1}{\theta}}(n)
\sum_{l=n+j}^\infty \widehat f(l)\widehat\theta(l-n-j).$$
Furthermore, \begin{align*}\sum_{n=0}^\infty \Bigl|\widehat{\frac{1}{\theta}}(n)\Bigr|
& \sum_{l=n+j}^\infty |\widehat f(l)||\widehat\theta(l-n-j)| \\
& \leq 
\sum_{n=0}^\infty \Bigl|\widehat{\frac{1}{\theta}}(n)\Bigr|
\Big(\sum_{l=n+j}^\infty |\widehat f(l)|^2\Big)^{1/2}
\Big(\sum_{l=n+j}^\infty |\widehat\theta(l-n-j)|^2\Big)^{1/2}  \\  
 & \leq
\sum_{n=0}^\infty \Bigl|\widehat{\frac{1}{\theta}}(n)\Bigr|\|S^{\ast n}f\|<\infty,\end{align*} 
and hence, 
$$\sum_{n=0}^\infty \widehat{\frac{1}{\theta}}(n)\sum_{l=n+j}^\infty \widehat f(l)\widehat\theta(l-n-j) = 
\sum_{l=j}^\infty\widehat f(l)\sum_{n=0}^{l-j}\widehat{\frac{1}{\theta}}(n)\widehat\theta(l-n-j)=
\widehat f(j)$$
(see \eqref{E:6.1}). We obtain that $(\theta(S^\ast)g, \chi^j)=\widehat f(j)$ for every $j\geq 0$,
that is, $\theta(S^\ast)g =f$. \end{proof}

\begin{corollary}\label{C:6.6} Suppose that $T\in\mathcal L(\mathcal H)$ is an a.c. polynomially bounded operator, $\theta$ is a singular inner function, 
$u_0\in\mathcal H$ is such that $$\sum_{n=0}^\infty \Bigl|\widehat{\frac{1}{\theta}}(n)\Bigr|\|T^{\ast n} u_0\|<\infty, \  \ 
\text{ and } \ \ 
u = \sum_{n=0}^\infty \widehat{\frac{1}{\theta}}(n)T^{\ast n} u_0.$$ Then $\theta(T^\ast)u=u_0$.\end{corollary}
\begin{proof} Since $\sum_{n=0}^\infty \bigl|\widehat{\frac{1}{\theta}}(n)\bigr|=\infty$, $\inf_n \|T^{\ast n} u_0\|=0$. 
 Since $T$ is power bounded, we conclude that $0=\inf_n \|T^{\ast n} u_0\|=\lim_n\|T^{\ast n} u_0\|$. 
Set $\mathcal H_0=\vee_{n\geq 0}T^{\ast n} u_0$. Then $\mathcal H_0\in\operatorname{Lat}T^\ast$. Set $T_0=T^\ast|_{\mathcal H_0}$.
Then $T_0$ is of class $C_{0\cdot}$. 

It is sufficient to prove that $\theta(T_0)u=u_0$. By \cite{bercpr}, there exists a contraction $R$ and a quasiaffinity $X$ such that 
$XT_0=RX$. It follows from this relation that $R$ is of class $C_{0\cdot}$. Therefore, there exists 
$\mathcal M\in \operatorname{Lat}S_\infty^\ast$ such that $R\cong S_\infty^\ast|_{\mathcal M}$  
{\cite[Theorem VI.2.3]{sznagy}}. 
Without loss of generality we can suppose that $R= S_\infty^\ast|_{\mathcal M}$. 
Put $v_0=Xu_0$. Then 
$\|R^n v_0\|=\|XT_0^nu_0\|\leq\|X\|\|T_0^nu_0\|$, therefore, 
$$\sum_{n=0}^\infty \Bigl|\widehat{\frac{1}{\theta}}(n)\Bigr|\|R^n v_0\|<\infty.$$
Set $$v=\sum_{n=0}^\infty \widehat{\frac{1}{\theta}}(n)R^n v_0.$$ By Lemma~\ref{L:6.5}, $\theta(R)v=v_0$. 
Also, $v=Xu$. We have \begin{align*}Xu_0=v_0= \theta(R)v = \theta(R)Xu=X\theta(T_0)u.\end{align*}Since $\ker X=\{0\}$, we conclude that $\theta(T_0)u=u_0$. 
\end{proof}

 \begin{theorem}\label{T:6.7} Suppose that  $T\in\mathcal L(\mathcal H)$ is an a.c. polynomially bounded operator, 
 $\tau\subset \mathbb T$ is a Borel set, $m(\tau)>0$, 
$X\in\mathcal L(\mathcal H, L^2(\tau,m))$, $XT=U(\tau)X$, and $\operatorname{clos}X\mathcal H=L^2(\tau,m)$.
Suppose  that $\theta$ is a singular inner function, $0\not\equiv g\in L^2(\tau,m)$,  and 
 \begin{align}\label{E:6.10}\sum_{n=0}^\infty \Bigl|\widehat{\frac{1}{\theta}}(n)\Bigr|\|T^{\ast n}X^\ast g \|<\infty. 
\end{align} 
Then the set $\{\xi\in\mathbb T :\ \ker\theta_\xi(T^\ast)\neq\{0\}\}$ contains a nonempty open subset of $\mathbb T$.
Thus, the operator $T$ has nontrivial hyperinvariant subspaces.
\end{theorem}
\begin{proof} For every $\xi\in\mathbb T$, define $u_\xi$ and $v_\xi$ as in Lemma~\ref{L:6.4}. By Corollary~\ref{C:6.6},  $\theta_\xi(T^\ast)u_\xi=X^\ast g$.
We have $$\theta_\xi(T^\ast)v_\xi=\theta_\xi(T^\ast)X^\ast (\theta_\xi)\widetilde{\ }g =
X^\ast  \theta_\xi(U(\tau)^\ast)(\theta_\xi)\widetilde{\ }g =X^\ast g.$$
Thus, $u_\xi-v_\xi\in\ker\theta_\xi(T^\ast)$ for every $\xi\in\mathbb T$. Let $\Xi$ be the set from Lemma~\ref{L:6.4}. Then
 $\{\xi\in\mathbb T :\ \ker\theta_\xi(T^\ast)\neq\{0\}\}\supset\mathbb T\setminus\Xi$. 
\end{proof}

 \begin{corollary}\label{C:6.8} Suppose  that $T\in\mathcal L(\mathcal H)$ is an a.c. polynomially bounded operator, 
 $\tau\subset \mathbb T$ is a Borel set, $m(\tau)>0$, 
$X\in\mathcal L(\mathcal H, L^2(\tau,m))$, $XT=U(\tau)X$, and $\operatorname{clos}X\mathcal H=L^2(\tau,m)$.
Suppose that  $0\not\equiv g\in L^2(\tau,m)$,   and
 $$\sum_{n=0}^\infty \|T^{\ast n}X^\ast g \|^2<\infty.$$ 
Then there exists a singular inner function $\theta$ such that $m(\operatorname{supp}\mu_\theta)=0$ and the set 
 $\{\xi\in\mathbb T :\ \ker\theta_\xi(T^\ast)\neq\{0\}\}$ contains a nonempty open subset of $\mathbb T$.\end{corollary}
\begin{proof} By Lemma~\ref{L:2.2}(iii) applied with $\varepsilon_n=\|T^{\ast n}X^\ast g \|$, there exists a dissymmetric weight $\omega$ such that 
$$\sum_{n=0}^\infty \|T^{\ast n}X^\ast g \|^2\omega(-n-1)^2<\infty.$$
By Theorem~\ref{T:2.1}, there exists a singular inner function $\theta$ such that 
\eqref{E:2.3} is fulfilled for $\theta$ and $\omega$ and $m(\operatorname{supp}\mu_\theta)=0$. 
We have \begin{align*}\sum_{n=0}^\infty  & \Bigl|\widehat{\frac{1}{\theta}}(n) \Bigr|\|T^{\ast n} X^\ast g\| \\  & \leq
\biggl(\sum_{n=0}^\infty \frac{1}{\omega(-n-1)^2}\Bigl|\widehat{\frac{1}{\theta}}(n)\Bigr|^2\biggr)^{1/2}
\biggl(\sum_{n=0}^\infty \|T^{\ast n}X^\ast g \|^2\omega(-n-1)^2\biggr)^{1/2}<\infty.\end{align*}
It remains to apply Theorem~\ref{T:6.7}.  \end{proof}

\begin{theorem}\label{T:6.9} Suppose  that $T\in\mathcal L(\mathcal H)$ is a power bounded operator, 
 $\tau\subset \mathbb T$ is a Borel set, $m(\tau)>0$, 
$X\in\mathcal L(\mathcal H, L^2(\tau,m))$, $XT=U(\tau)X$, and $\operatorname{clos}X\mathcal H=L^2(\tau,m)$.
Suppose  that $\theta$ is a singular inner function, $0\not\equiv\varphi\in A^+(\mathbb T)$,  
$\theta\varphi\in A^+(\mathbb T)$, $0\not\equiv g\in L^2(\tau,m)$, and 
the condition \eqref{E:6.10} is fulfilled for $\theta$ and $g$. 
Then the set $\{\xi\in\mathbb T :\ \ker(\theta\varphi)_\xi(T^\ast)\neq\{0\}\}$ contains a  nonempty open subset of $\mathbb T$. Thus, the operator $T$ has nontrivial hyperinvariant subspaces.\end{theorem}
\begin{proof}
For every $\xi\in\mathbb T$, define $u_\xi$ and $v_\xi$ as in Lemma~\ref{L:6.4}. 
Clearly, $\varphi_\xi$, $(\theta\varphi)_\xi\in A^+(\mathbb T)$ for every
  $\xi\in\mathbb T$.
Set $M=\sup_{k\geq 0}\|T^k\|$. We have 
\begin{equation*}\begin{gathered}\sum_{k=0}^\infty\sum_{n=0}^\infty
\Bigl\|\widehat{(\theta\varphi)}_\xi(k)\widehat{\frac{1}{\theta}}(n)T^{\ast (k+n)}X^\ast g \Bigr\|\leq
\sum_{k=0}^\infty\sum_{n=0}^\infty
|\widehat{(\theta\varphi)}_\xi(k)|\Bigl|\widehat{\frac{1}{\theta}}(n)\Bigr|M\|T^{\ast n}X^\ast g \| \\
=M\sum_{k=0}^\infty|\widehat{(\theta\varphi)}_\xi(k)|
\sum_{n=0}^\infty\Bigl|\widehat{\frac{1}{\theta}}(n)\Bigr|\|T^{\ast n}X^\ast g \|<\infty, \end{gathered}\end{equation*}
therefore, \begin{equation*}\begin{gathered}
(\theta\varphi)_\xi(T^\ast)u_\xi=
\sum_{k=0}^\infty\widehat{(\theta\varphi)}_\xi(k)T^{\ast k}
\sum_{n=0}^\infty \widehat{\frac{1}{\theta_\xi}}(n)T^{\ast n}X^\ast g \\ 
=\sum_{n=0}^\infty\Big(\sum_{k=0}^n
\widehat{(\theta\varphi)}_\xi(n-k) \widehat{\frac{1}{\theta_\xi}}(k)\Big)T^{\ast n}X^\ast g
=\sum_{n=0}^\infty\widehat\varphi_\xi(n)T^{\ast n}X^\ast g=\varphi_\xi(T^\ast)X^\ast g. \end{gathered}\end{equation*}
Furthermore, \begin{align*}(\theta\varphi)_\xi(T^\ast)v_\xi=
(\theta\varphi)_\xi(T^\ast)X^\ast (\theta_\xi)\widetilde{\ } g
&=X^\ast(\theta\varphi)_\xi(U(\tau)^\ast)(\theta_\xi)\widetilde{\ } g \\
 = X^\ast\varphi_\xi(U(\tau)^\ast)\theta_\xi(U(\tau)^\ast)(\theta_\xi)\widetilde{\ } g&=
X^\ast\varphi_\xi(U(\tau)^\ast) g=\varphi_\xi(T^\ast) X^\ast g.\end{align*}
We obtain that  $u_\xi-v_\xi\in\ker(\theta\varphi)_\xi(T^\ast)$ for every $\xi\in\mathbb T$. 
 Let $\Xi$ be the set from Lemma~\ref{L:6.4}. Then
 $\{\xi\in\mathbb T :\ \ker\theta_\xi(T^\ast)\neq\{0\}\}\supset\mathbb T\setminus\Xi$. \end{proof}

\begin{corollary}\label{C:6.10} Suppose  that $T\in\mathcal L(\mathcal H)$ is a power bounded operator, 
 $\tau\subset \mathbb T$ is a Borel set, $m(\tau)>0$, 
$X\in\mathcal L(\mathcal H, L^2(\tau,m))$, $XT=U(\tau)X$, and $\operatorname{clos}X\mathcal H=L^2(\tau,m)$.
Suppose  that there exist $0\not\equiv g\in L^2(\tau,m)$,   a sequence $\{w_n\}_{n\geq 1}$ of positive numbers, and $C>0$  such that 
$\{w_n\}_{n\geq 1}$ satisfies  the conditions of Theorem~\ref{T:2.4}, and  
$$\|T^{\ast n}X^\ast g \|\leq C/ w_{n+1} \  \  \text{ for sufficiently large } n.$$
Then there exist  a singular inner function  $\theta$ and an outer function $\varphi\in A^+(\mathbb T)$ such that 
$\operatorname{supp}\mu_\theta$ is a Carleson set, $m(\operatorname{supp}\mu_\theta)=0$,
$\theta\varphi\in A^+(\mathbb T)$,  and  the set $\{\xi\in\mathbb T :\ \ker(\theta\varphi)_\xi(T^\ast)\neq\{0\}\}$ 
contains a nonempty open subset of $\mathbb T$.\end{corollary}
\begin{proof} Let $0<s<1$. By Theorem~\ref{T:2.4} and Lemma~\ref{L:2.5}(ii), there exist   a singular inner function  $\theta$  and $C_1>0$ such that 
$\operatorname{supp}\mu_\theta$ is a Carleson set, $m(\operatorname{supp}\mu_\theta)=0$, and 
$$ \Bigl|\widehat{\frac{1}{\theta}}(n-1)\Bigr|\leq C_1 w_n^s \ \ \text{ for all } n\geq 1.$$
We have $$\sum_{n=0}^\infty  \Bigl|\widehat{\frac{1}{\theta}}(n)\Bigr|\|T^{\ast n}X^\ast g \|\leq
CC_1\sum_{n=0}^\infty\frac{1}{ w_{n+1}^{1-s}}<\infty$$
by  Lemma~\ref{L:2.5}(ii). Thus, condition \eqref{E:6.10} is fulfilled for $\theta$. 
Since $\operatorname{supp}\mu_\theta$ is a Carleson set and $m(\operatorname{supp}\mu_\theta)=0$, there exists 
an outer function $\varphi\in A^+(\mathbb T)$ such that $\theta\varphi\in A^+(\mathbb T)$ \cite{tawi}. It remains to apply Theorem~\ref{T:6.9}.  
\end{proof}

Theorems~\ref{T:6.7} and \ref{T:6.9} and Corollaries~\ref{C:6.8} and~\ref{C:6.10} are formulated for an arbitrary operator $X$ (with dense range) which intertwines 
$T$ with $U(\tau)$. The following proposition shows how Theorems~\ref{T:6.7} and \ref{T:6.9} and Corollaries~\ref{C:6.8} and~\ref{C:6.10} can be formulated 
in terms of the unitary asymptote of $T$. 

\begin{proposition}\label{P:6.11} Suppose that  $T\in\mathcal L(\mathcal H)$ is an operator,  $U\in\mathcal L(\mathcal K)$ is a unitary operator, $(Y,U)$ is the unitary asymptote 
of $T$, and $\operatorname{clos}Y\mathcal H=\mathcal K$. If $U$ is a.c., then for every $0\neq y\in\mathcal K$ there exist  a Borel set 
$\tau\subset \mathbb T$, $X\in\mathcal L(\mathcal H, L^2(\tau,m))$, and $0\not\equiv g\in L^2(\tau,m)$ such that $m(\tau)>0$, 
 $XT=U(\tau)X$,  $\operatorname{clos}X\mathcal H=L^2(\tau,m)$, and $X^\ast g = Y^\ast y$. 
Conversely, if there exist  a Borel set $\tau\subset \mathbb T$   and $X\in\mathcal L(\mathcal H, L^2(\tau,m))$    such that $m(\tau)>0$, 
$XT=U(\tau)X$, and $\operatorname{clos}X\mathcal H=L^2(\tau,m)$, then for every $0\not\equiv g\in L^2(\tau,m)$, 
there exists $0\neq y\in\mathcal K$ such that $X^\ast g = Y^\ast y$. \end{proposition}
\begin{proof} For $0\neq y\in\mathcal K$ put $\mathcal M=\vee_{n\in\mathbb Z}U^ny$. Then $\mathcal M$ is a reducing subspace of $U$, and 
$\operatorname{clos}P_{\mathcal M}Y\mathcal H=\mathcal M$. If $U$ is a.c., there exist  a Borel set $\tau\subset \mathbb T$ such that $m(\tau)>0$ and 
a unitary operator $W\in\mathcal L(\mathcal M, L^2(\tau,m))$ such that $WU|_{\mathcal M}=U(\tau)W$. Put $g=Wy$ and $X=WP_{\mathcal M}Y$. Then 
$g\not\equiv 0$, the equality $XT= U(\tau)X$ follows from the intertwining properties of $W$ and $Y$ and the reducing property of $\mathcal M$, 
and $X^\ast g = Y^\ast y$ because $W^\ast=W^{-1}$.  

Conversely, suppose that there exist  a Borel set $\tau\subset \mathbb T$  and 
$X$ such that $m(\tau)>0$ and $XT=U(\tau)X$. By the definition of the unitary asymptote, there exists an operator 
$Z$ such that $X=ZY$ and $ZU=U(\tau)Z$. For $g\in L^2(\tau,m)$, put $y=Z^\ast g$. Then $X^\ast g = Y^\ast y$.
If  $\operatorname{clos}X\mathcal H=L^2(\tau,m)$, then 
$\operatorname{clos}Z\mathcal K=L^2(\tau,m)$, therefore, $\ker Z^\ast=\{0\}$. Thus, if $g\not\equiv 0$, then $y\neq 0$. \end{proof}

\begin{remark}\label{R:6.12} In Sec. 3 and 4 of the present paper, we use the functional calculus for power bounded operators on the algebra 
$A^+(\mathbb T)$ of analytic functions with absolutely summable Taylor coefficients.
In \cite{pel}, a  functional calculus for power bounded operators on some larger  algebras of analytic functions 
is constructed. However, the author does not have enough knowledge on these algebras  to extend to them the arguments used for $A^+(\mathbb T)$.\end{remark}

\section{Example of a quasianalytic contraction}

In this section, an example of a quasianalytic contraction $T$ is constructed such that $\sigma(T)=\mathbb T$ and 
$\emptyset\neq\pi(T)\neq \mathbb T$. Actually, an operator similar to a contraction is constructed, but it is easy to see 
that the residual and quasianalytic spectral sets of operators having unitary asymptote are preserved under similarity. 

\subsection{Polynomially bounded operators: quasianalyticity and  full analytic invariant subspaces}

Before  constructing an example, we formulate Theorem~\ref{T:4.6} (which is actually from \cite{rej}), which shows that the constructed example will be typical.
 Namely, to construct a needed quasianalytic operator $T$ we first  take an appropriate quasianalytic operator $T_1\in\mathcal H_1$ such that $\mathcal H_1$ is a full analytic invariant subspace for $T_1$ (see Definition ~\ref{D:4.5}); therefore, $\sigma(T_1)=\operatorname{clos}\mathbb D$. The operator $T_1$ will be the restriction of the needed operator $T$ on its invariant subspace. 

The following lemma is well known, see, for example, {\cite[Lemma 3.1]{berc}} or {\cite[Lemma 6]{tak}}. Detailed proof can be found in {\cite[Lemma 2.2]{gaminner}}.

\begin{lemma} \label{L:4.2} Suppose   that $C>0$, $k\in\mathbb N$,  and an operator $T$ 
is such that $\sigma(T)\subset\operatorname{clos}\mathbb D$ and $\|(T-\lambda I)^{-1}\|\leq C/(|\lambda|-1)^k$ for all 
$\lambda\notin\operatorname{clos}\mathbb D$. 
Put $$\Lambda =\sigma(T)\cup\bigl\{\lambda\in\mathbb D: \ \ \|(T-\lambda I)^{-1}\|\geq C/(1-|\lambda|)^k\bigr\}$$
and $$\mathcal Z=\bigl\{\zeta\in\sigma(T)\cap\mathbb T: \ \ \sup\{r\in [0,1): r\zeta\in\Lambda\}<1\bigr\}.$$ 
If $\mathcal Z$ is uncountable, then there exist $\mathcal M_1$, $\mathcal M_2\in\operatorname{Hlat}T$ such that $\mathcal M_1\neq\{0\}$, 
$\mathcal M_2\neq\{0\}$, and $\sigma(T|_{\mathcal M_1})\cap\sigma(T|_{\mathcal M_2})=\emptyset$.\end{lemma}

Recall the definition.

\begin{definition}\label{D:4.5} Let $T$ be an operator (on a Hilbert space), and let $\mathcal M\in\operatorname{Lat}T$. The space $\mathcal M$ is called 
a {\it full analytic invariant subspace} of $T$ if there exists a nonzero analytic function  $\mathbb D\to \mathcal M$, 
$\lambda\mapsto e_\lambda$, 
such that 
$$(T|_{\mathcal M}- \overline\lambda I_{\mathcal M})^\ast e_\lambda =0 \ \ \text { for every } \lambda\in\mathbb D, \ \ \text{ and } \
\bigvee_{\lambda\in\mathbb D} { e_\lambda}=\mathcal M.$$\end{definition}

The following theorem is a straightforward consequence of the results of \cite{rej}.

\begin{theorem} \label{T:4.6}Suppose  that $T\in\mathcal L(\mathcal H)$ is an a.c. polynomially bounded operator such that $\mathbb T\subset\sigma(T)$ and $\pi(T)\neq\emptyset$. 
Then the set of elements of $\mathcal H$ which generate a full analytic invariant subspace for $T$ is dense in $\mathcal H$, 
and $T$ is reflexive.\end{theorem}
\begin{proof}  To apply results of \cite{rej} we need to show that 
$T$ is of class $C_{\cdot 0}$, and  there exists $C'>0$ such that
for all $C>C'$ the set of points in $\mathbb T$ which are not radial limits of points in $\Lambda_C$ is at most countable, where \begin{align*}
\Lambda_C =\Bigl\{\lambda\in\mathbb D: \ &\text{ there exists } x\in\mathcal H \ \text{ such that } \|x\|=1 \\
 & \text{ and } \ \  \|(T^\ast-\lambda I)x\|<\frac{(1-|\lambda|)^2}{C}\Bigr\}.\end{align*}
By {\cite[Proposition 33]{ker5}}, $T$ is of class $C_{10}$. In particular, $T$ has no eigenvalues. Therefore, 
$$\Lambda_C =\bigl\{\lambda\in \mathbb D :\ \ \overline\lambda \in\sigma(T) \ \ \text{ or } \ \ \|(T-\overline\lambda I)^{-1}\|> C/(1-|\lambda|)^2\bigr\}.$$
If we suppose that the set of points in $\mathbb T$ which are not radial limits of points in $\Lambda_C$ is uncountable for some $C>0$,
we can apply Lemma~\ref{L:4.2}   with some $C_1>C$ to $T$. We obtain $\mathcal M_1$, $\mathcal M_2\in\operatorname{Hlat}T$ such that $\mathcal M_1\neq\{0\}$, 
$\mathcal M_2\neq\{0\}$, and $\sigma(T|_{\mathcal M_1})\cap\sigma(T|_{\mathcal M_2})=\emptyset$. 
By {\cite[Proposition 35]{ker5}}, $\pi(T)\subset\sigma(T|_{\mathcal M_l})$  for $l=1,2$. Consequently,  
$\pi(T)\subset\sigma(T|_{\mathcal M_1})\cap\sigma(T|_{\mathcal M_2})=\emptyset$, a contradiction. Therefore, $T$ satisfies 
the conditions of {\cite[Theorem 6.5]{rej}}, and the conclusion of the theorem follows from the proof of 
{\cite[Theorem 6.5]{rej}}. \end{proof}

\subsection{Similarity to a contraction}

For a natural number $N$ a $N\times N$ matrix can be regarded as an operator on 
$\ell_N^2$, its norm is denoted by the symbol $\|\cdot\|_{\mathcal L(\ell_N^2)}$. 
For a family of polynomials $[\varphi_{ij}]_{i,j=1}^N$ put
$$\| [\varphi_{ij}]_{i,j=1}^N\|_{H^\infty(\ell_N^2)}= 
\sup\bigl\{\| [\varphi_{ij}(z)]_{i,j=1}^N\|_{\mathcal L(\ell_N^2)} \ : \ z\in\mathbb D\bigr\}.$$
For an operator  $T \in\mathcal L(\mathcal H)$  and   a  family   of  polynomials 
$[\varphi_{ij}]_{i,j=1}^N$  
the operator  $$[\varphi_{ij}(T)]_{i,j=1}^N\in\mathcal L\bigl(\oplus_{j=1}^N \mathcal H\bigr)$$ 
is defined.
$T$ is called {\it completely polynomially bounded}, 
if there exists a constant $M$ such that
\begin{align}\label{E:7.1}\| [\varphi_{ij}(T)]_{i,j=1}^N\|\leq M\| [\varphi_{ij}]_{i,j=1}^N\|_{H^\infty(\ell_N^2)}
\end{align}
for every family of polynomials  $[\varphi_{ij}]_{i,j=1}^N$ and every $N\geq 1$. 
If $T$ is a contraction, then $M=1$ in \eqref{E:7.1}.
The following criterion for an operator to be similar to a contraction is proved in \cite{pau}.

{\it An operator $T$ is similar to a contraction if and only if $T$ is completely polynomially bounded. }

For an index $k\geq 0$ and a polynomial $\varphi$ put 
\begin{align}\label{E:7.2}(\varphi)_k(z)=\sum_{n\geq k+1}\widehat\varphi(n)z^{n-k-1}, \ \ z\in\mathbb C. \end{align}
The following lemma is a particular case of {\cite[Lemma 3.1]{pruvost}}.

\begin{lemma}[\cite{pruvost}] \label{L:7.1} There is a constant $C>0$ such that for every  $N\geq 1$, $k\geq 0$ and for every family of polynomials $[\varphi_{ij}]_{i,j=1}^N$, 
$$\| [(\varphi_{ij})_k]_{i,j=1}^N\|_{H^\infty(\ell_N^2)}\leq C\log (k+2) \| [\varphi_{ij}]_{i,j=1}^N\|_{H^\infty(\ell_N^2)}.$$
\end{lemma}

Recall that the definition, notations and references on (well-known and intensively studied) weighted shifts were given in Sec. 2.

\begin{theorem}\label{T:7.2} Suppose that $T_1\in \mathcal L(\mathcal H_1)$ is a contraction, $T_1$ has no eigenvalues, and there exists $0\neq x_0\in \mathcal H_1$ such that 
$\mathcal H_1=(T_1-\lambda I)\mathcal H_1\dotplus \mathbb C x_0$ for every $\lambda\in \mathbb D$ 
(the norm of the (nonorthogonal) projection defined by $\dotplus$ may depend on $\lambda$). Suppose  that $\omega$ is 
a submultiplicative dissymmetric weight, 
and \begin{align}\label{E:7.8}\sum_{n=1}^\infty \Bigl( \frac{\log n}{\omega(-n)}\Bigr) ^2 <\infty.\end{align}
Define $A\in\mathcal L(\ell^2_{\omega -}, \mathcal H_1)$ by the formula $Au = u(-1)x_0$, $u\in\ell^2_{\omega -}$.
 Put $$T=\begin{pmatrix} T_1 & A \cr \mathbb O & S_{\omega -}\end{pmatrix}, \ \ \  \ T\in\mathcal L(\mathcal H_1\oplus\ell^2_{\omega -}).$$
Then $\sigma(T)=\mathbb T$ and $T$ is similar to a contraction.\end{theorem}
\begin{proof} For a polynomial $\varphi$, set $A_\varphi = P_{\mathcal H_1\oplus \{0\}}\varphi(T)|_{\{0\}\oplus \ell^2_{\omega -}}$, 
that is, $A_\varphi$ is the operator from the right upper corner of the matrix form of $\varphi(T)$. It is easy to see that 
\begin{align}\label{E:7.9}A_\varphi u= \sum_{k=0}^\infty u(-1-k)(\varphi)_k(T_1)x_0, \ \ \ u\in\ell^2_{\omega -}, 
\end{align}
where $(\varphi)_k$ is defined by \eqref{E:7.2} (since  $\varphi$ is a polynomial, the sum in \eqref{E:7.9} is actually finite).

Since $T_1$ and $S_{\omega -}$ are contractions, 
to show that $T$ satisfies \eqref{E:7.1} it is sufficient to show that there exists a constant $C_1>0$ such that 
\begin{align}\label{E:7.10}\|[A_{\varphi_{ij}}]_{i,j=1}^N [u_j]_{j=1}^N\|\leq C_1\|[\varphi_{ij}]_{i,j=1}^N\|_{H^\infty(\ell^2_N)}\|[u_j]_{j=1}^N\| \end{align}
for every family $[u_j]_{j=1}^N$ of elements of $\ell^2_{\omega -}$, for every family $[\varphi_{ij}]_{i,j=1}^N$ of 
polynomials, and every index $N\geq 1$. 
Let  $[u_j]_{j=1}^N$ be a family of elements of $\ell^2_{\omega -}$, put 
$$a_{jk}=u_j(-1-k)\omega(-1-k), \ \ \ k\geq 0,  \ \ 1\leq j\leq N.$$ 
Then \begin{align}\label{E:7.11}\|[u_j]_{j=1}^N\|^2=\sum_{k= 0}^\infty\sum_{j=1}^N|a_{jk}|^2. \end{align}
We infer from \eqref{E:7.9} that 
\begin{align*}[A_{\varphi_{ij}}]_{i,j=1}^N [u_j]_{j=1}^N& =\Bigl [\sum_{j=1}^N A_{\varphi_{ij}}u_j\Bigr]_{i=1}^N = 
\Bigl[\sum_{j=1}^N \sum_{k=0}^\infty u_j(-1-k)(\varphi_{ij})_k(T_1)x_0\Bigr]_{i=1}^N  \\
& = \sum_{k=0}^\infty\frac{1}{\omega(-1-k)}\Bigl[\sum_{j=1}^N a_{jk}(\varphi_{ij})_k(T_1)x_0\Bigr]_{i=1}^N \\
& = \sum_{k=0}^\infty\frac{1}{\omega(-1-k)}[(\varphi_{ij})_k(T_1)]_{i,j=1}^N [a_{jk}x_0]_{j=1}^N.\end{align*}
Applying Lemma~\ref{L:7.1} and \eqref{E:7.1} for the contraction $T_1$, we obtain that 
\begin{align*} & \|[A_{\varphi_{ij}}]_{i,j=1}^N  [u_j]_{j=1}^N\|  \leq 
\sum_{k=0}^\infty\frac{1}{\omega(-1-k)}\|[(\varphi_{ij})_k(T_1)]_{i,j=1}^N\|\|[a_{jk}x_0]_{j=1}^N\| \\
&
\leq \sum_{k=0}^\infty\frac{1}{\omega(-1-k)}\|[(\varphi_{ij})_k]_{i,j=1}^N\|_{H^\infty(\ell^2_N)}
\biggl(\sum_{j=1}^N|a_{jk}|^2\biggr)^{1/2}\|x_0\| \\
& \leq \sum_{k=0}^\infty\frac{1}{\omega(-1-k)}C\log (k+2)\|[\varphi_{ij}]_{i,j=1}^N\|_{H^\infty(\ell^2_N)}
\biggl(\sum_{j=1}^N|a_{jk}|^2\biggr)^{1/2}\|x_0\| \\
& \leq C\|[\varphi_{ij}]_{i,j=1}^N\|_{H^\infty(\ell^2_N)}\|x_0\|
\biggl(\sum_{k=0}^\infty\Bigl(\frac{\log (k+2)}{\omega(-1-k)}\Bigr)^2\biggr)^{1/2}
\biggl(\sum_{k=0}^\infty\sum_{j=1}^N|a_{jk}|^2\biggr)^{1/2}.\end{align*}
Now estimate \eqref{E:7.10} follows from the latter estimate and \eqref{E:7.11}.

Since $T$ is similar to a contraction,  $\sigma(T)\subset\operatorname{clos}\mathbb D$. Let $\lambda\in\mathbb D$, 
and let $x\in \mathcal H_1$ and $u\in\ell^2_{\omega -}$ be such that $(T-\lambda I)(x\oplus u)=0$. Then 
$S_{\omega -}u = \lambda u$, therefore, $u(n)=\lambda^{-n-1}u(-1)$ for $n\leq -1$. Furthermore, 
$u(-1)x_0=-(T_1-\lambda I)x$. By the condition on $T_1$, $u(-1)=0$, therefore, $x=0$ and $u=0$. Thus, $T$ has no eigenvalues in $\mathbb D$. 

By {\cite[Proposition 2.3]{est}}, $\sigma(S_\omega)\subset\mathbb T$. Therefore, $(S_\omega-\lambda I)\ell^2_\omega = \ell^2_\omega$ 
for every  $\lambda\in\mathbb D$. Taking into account 
the definitions of $T$,  $S_\omega$, and the assumption on $T_1$, we obtain that 
$(T-\lambda I)(\mathcal H_1\oplus\ell^2_{\omega -})=\mathcal H_1\oplus\ell^2_{\omega -}$ for every  $\lambda\in\mathbb D$. 
Thus, $\sigma(T)\subset\mathbb T$. 

By {\cite[Theorem 0.8]{rara}}, the spectrum of the restriction of an operator on its invariant subspace is contained in the polynomially 
convex hull of  its spectrum. Therefore, if $\sigma(T)\neq \mathbb T$, then  $\sigma(T_1)\neq\operatorname{clos} \mathbb D$, 
a contradiction with the assumption on $T_1$. Thus,  $\sigma(T)=\mathbb T$. \end{proof}

\subsection{Subnormal operators}

Let $\nu$ be a positive finite Borel measure on $\operatorname{clos}\mathbb D$. Denote by 
$P^2(\nu)$ the closure of analytic polynomials in $L^2(\nu)$, and by $S_\nu$
the operator of multiplication by the independent variable in $P^2(\nu)$, i.e.
$$S_\nu\in\mathcal L(P^2(\nu)), \ \ \ \ (S_\nu f)(z)=zf(z), \ \ f\in P^2(\nu),
 \ \ z\in \operatorname{clos}\mathbb D.$$
Clearly, $S_\nu$ is a contraction. 

Denote   by $m_2$ the normalized planar Lebesgue measure on $\mathbb D$. 
For $-1<\alpha\leq 0$ put 
$$ \text{\rm d}m_{2,\alpha}(z)=(1-|z|^2)^\alpha\text{\rm d}m_2(z) $$
and \begin{align}\label{E:7.12} v_\alpha\colon\mathbb Z_+\to (0,\infty), \ \ \ \ v_\alpha(n)^2=\frac{1}{(n+1)^{\alpha+1}}, \ \ \ n\in\mathbb Z_+. \end{align}
The space $P^2(m_{2,\alpha})$ is the {\it weighted Bergman space}. It is well-known that  
$P^2(m_{2,\alpha})$ is the space of all functions in $L^2(m_{2,\alpha})$ that are analytic in $\mathbb D$, 
and  \begin{align} \label{E:7.13} \|f\|^2_{P^2(m_{2,\alpha})}\asymp\sum_{n=0}^\infty|\widehat f(n)|^2 v_\alpha(n)^2,  \ \ \ f\in P^2(m_{2,\alpha}) \end{align}
(the estimate in \eqref{E:7.13} depends on $\alpha$), see, for example, {\cite[Ch.1.1,1.2]{heko}}. In other words, the operator 
 \begin{align} \label{E:7.14} J_\alpha\in\mathcal L(P^2(m_{2,\alpha}), \ell^2_{v_\alpha +}), \ \ \ J_\alpha f = \{\widehat f(n)\}_{n=0}^\infty,  
\ \ f\in P^2(m_{2,\alpha}), \end{align}
is invertible, and, clearly, $J_\alpha S_{m_{2,\alpha}} = S_{v_\alpha +}J_\alpha $. 

The following theorem is proved in {\cite[Sec. 2]{gama}} for $\alpha=0$, but the proof is the same for all $-1<\alpha\leq 0$. 
We refer to this paper because its results are formulated in a form convenient for the purpose of the present paper. The reader 
interested in subnormal operators can see \cite{conw2}, \cite{alers} and references therein.    

 \begin{theorem}\label{T:7.3} Let $-1<\alpha\leq 0$. Suppose that $\tau\subset\mathbb T$ is a Borel set such that $0<m(\tau)<1$. Then there exists a positive 
finite Borel measure $\mu$ on $\operatorname{clos}\mathbb D$ such that the space $P^2(m_{2,\alpha}+\mu)$ and the operator 
$S_{m_{2,\alpha}+\mu}$ have the following properties:

\begin{enumerate}[\upshape (i)]
\item the mapping $$f\mapsto f(\lambda), \ \ P^2(m_{2,\alpha}+\mu)\to \mathbb C$$ is (linear, bounded) functional on 
$P^2(m_{2,\alpha}+\mu)$ for every $\lambda\in\mathbb D$, every function $f\in P^2(m_{2,\alpha}+\mu)$ is analytic in $\mathbb D$, 
and has nontangential boundary values $f(\zeta)$ for almost all $\zeta\in\tau$ with respect to $m$;

\item the contraction $S_{m_{2,\alpha}+\mu}$ is of class $C_{10}$, and $(Y_\tau, U(\tau))$ is its unitary asymptote,
where $Y_\tau\in\mathcal L(P^2(m_{2,\alpha}+\mu),L^2(\tau,m))$ acts by the formula $(Y_\tau f)(\zeta)= f(\zeta)$ for a.e. $\zeta\in\tau$;

\item $P^2(m_{2,\alpha}+\mu) = (S_{m_{2,\alpha}+\mu}-\lambda I)P^2(m_{2,\alpha}+\mu) \dotplus \mathbb C\text{{\bf 1}}$ 
for every $\lambda\in\mathbb D$, where $\text{{\bf 1}}\in P^2(m_{2,\alpha}+\mu) $, $\text{{\bf 1}}(z)=1$ for all
 $z\in\operatorname{clos}\mathbb D$.\end{enumerate}\end{theorem}  
 \begin{proof} The conclusions (i) and (ii) of the theorem are proved in [Gam, Sec. 2]. 
Since functions from $P^2(m_{2,\alpha}+\mu)$ are analytic, \begin{align}\label{E:7.15}\text{{\bf 1}}\notin (S_{m_{2,\alpha}+\mu}-\lambda I)P^2(m_{2,\alpha}+\mu)   \end{align} 
for every $\lambda\in\mathbb D$. It is  proved in {\cite[Sec. 2]{gama}} that $I-S_{m_{2,\alpha}+\mu}^\ast S_{m_{2,\alpha}+\mu}$ is compact. 
Therefore,  $S_{m_{2,\alpha}+\mu}$  is  a  compact  perturbation  of   an  isometry,   consequently, 
 $S_{m_{2,\alpha}+\mu}-\lambda I$ is bounded below for every $\lambda\in\mathbb D$. Thus,
$(S_{m_{2,\alpha}+\mu}-\lambda I)P^2(m_{2,\alpha}+\mu)$ is closed. Since $S_{m_{2,\alpha}+\mu}$ is cyclic,
 $\operatorname{codim}(S_{m_{2,\alpha}+\mu}-\lambda I)\leq 1$ for every $\lambda\in\mathbb D$. Now (iii) follows from \eqref{E:7.15}.

Also (iii) can be proved using that for every compact set  $K\subset \mathbb D$ there exists a constant $C>0$ 
(which depends on $K$)  such that 
$$\|f\|_{P^2(m_{2,\alpha}+\mu)}\leq C\|f\|_{L^2(\operatorname{clos}\mathbb D\setminus K,m_{2,\alpha}+\mu)}$$ for every function 
$f\in P^2(m_{2,\alpha}+\mu)$, because functions from $P^2(m_{2,\alpha}+\mu)$ are analytic, see \cite{alers}.   \end{proof}

\subsection{Properties of the constructed operator}

In this subsection we finish the construction of an example and studies of its properties. 

  \begin{proposition}\label{P:7.4} Suppose that $\tau\subset\mathbb T$ is a Borel set such that $0<m(\tau)<1$, $\omega$ is 
a submultiplicative dissymmetric weight satisfying \eqref{E:7.8}, and $T$ is the operator from Theorem~\ref{T:7.2} with $T_1=S_{m_{2,\alpha}+\mu}$ 
and $x_0=\text{{\bf 1}}$, where $S_{m_{2,\alpha}+\mu}$ and $\text{{\bf 1}}$ 
are from  Theorem~\ref{T:7.3}. Let $\mathcal J_{\omega -}\in \mathcal L(\ell^2_{\omega -}, H^2_-)$  be the restriction of the natural imbedding from \eqref{E:2.2} 
on $\ell^2_{\omega -}$. Define $X$ as follows:
\begin{align*} X\in\mathcal L\bigl(P^2(m_{2,\alpha}+\mu)\oplus\ell^2_{\omega -}, L^2(\tau,m)\bigr), & \ \ \ X(f\oplus u) =Y_\tau f + (\mathcal J_{\omega -}u)|_\tau, \\
& f\in P^2(m_{2,\alpha}+\mu),   \ \ u\in \ell^2_{\omega -}. \end{align*}
Then $(X, U(\tau))$ is a unitary asymptote of $T$.\end{proposition}
\begin{proof}  Since $T$ is similar to a contraction, $T$ has a unitary asymptote; denote it by $(Y, U)$. 
Recall that $YT=UY$. Since $S_{\omega -}$ is of class $C_{0\cdot}$, by \cite{ker1} and  Theorem~\ref{T:7.3}(ii), $U\cong U(\tau)$. Without loss of generality, suppose  that 
$U = U(\tau)$. By Theorem~\ref{T:7.3}(ii) and 
\cite{ker1}, $Y|_{P^2(m_{2,\alpha}+\mu)} = Y_\tau$. For $n\leq -1$ define $u_n\in\ell^2_{\omega -}$ by the formula
 $u_n(n)=1$, $u_n(k)=0$ for $k\leq -1$, $k\neq n$. It is easy to see from \eqref{E:7.9} that $T^{-n}(0\oplus u_n)=\text{{\bf 1}}\oplus 0$.
Therefore, $$\chi^0|_\tau = Y_\tau\text{{\bf 1}} = Y(\text{{\bf 1}}\oplus 0) = YT^{-n}(0\oplus u_n)$$ 
(where $\chi(\zeta)=\zeta$, $\zeta\in\mathbb T$). Furthermore, since $T$ is invertible, we have 
$$ Y(0\oplus u_n) = YT^nT^{-n}(0\oplus u_n) = U(\tau)^nYT^{-n}(0\oplus u_n) =  U(\tau)^n(\chi^0|_\tau) = \chi^n|_\tau$$ 
for $n\leq -1$.  Thus, $Y=X$.\end{proof}

The following theorem  is a straightforward consequence of \cite{kell3}. It gives sufficient conditions on the operator from 
Proposition~\ref{P:7.4} to  have nontrivial hyperinvariant subspaces. The author does not know whether  the operator from 
Proposition~\ref{P:7.4} satisfies  the conditions of Theorem~\ref{T:6.7}.   

 \begin{theorem}\label{T:7.5} Let $\omega$  and $T$ be as in Proposition~\ref{P:7.4}. Let  $\omega$ satisfy 
one of the following assumptions.
\begin{enumerate}[\upshape (i)]
\item  There exist $C>0$ and a nonincreasing sequence $\{w(n)\}_{n\leq -1}$  such that 
\begin{align}\label{E:7.16} w(n)\leq C\omega(n) \ \ \text{ for all }  n\leq -1, \end{align}
and $w$ has the following properties: 
$w(n)\geq 1$ for all $n\leq -1$, \begin{align*}\lim_{n\to\infty} w(-n)^{1/n} =1, \ \ 
w(-n)/n^c\to\infty \text{ when } n\to +\infty, \end{align*} the sequences $\Bigl\{\frac{w(-n-1)}{w(-n)}\Bigr\}_{n=1}^\infty$, $\Bigl\{\frac{w(-n-1)}{w(-n)}\frac{n^c}{(n+1)^c}\Bigr\}_{n=1}^\infty$,  and $\Bigl\{\frac{\log w(-n)}{n^b}\Bigr\}_{n=1}^\infty$
are nonincreasing for some $c>1/2$  and $b<1/2$, and $$\sum_{n\geq 1}\frac{1}{n \log w(-n)}<\infty;$$ 

\item  $\log\omega(-n)/\sqrt n\to \infty$ when $n\to +\infty$.
\end{enumerate}
Then there exists a singular inner function $\theta$ such that $\operatorname{ran}\theta(T)$ is not dense.\end{theorem}
\begin{proof} Set $\omega'(n)=\omega(n)$,  $w'(n)=w(n)$, $n\leq -1$, $w'(n)=\omega'(n)=v_\alpha(n)$, $n\geq 0$, where $v_\alpha$ is 
defined in \eqref{E:7.12}. Then $w'$ and $\omega'$ are nonincreasing functions on $\mathbb Z$. Let 
$$ J_{\alpha,\mu}\in\mathcal L\bigl(P^2(m_{2,\alpha}+\mu), P^2(m_{2,\alpha})\bigr)$$
be a natural imbedding. Recall that $J_\alpha$ is defined in \eqref{E:7.14}. 
It easy to see that the operator 
$J_\alpha J_{\alpha,\mu}\oplus I_{\ell^2_{\omega -}}$ realizes the relation 
$T\prec S_{\omega'}$. 

If  condition (ii) is fulfilled, then, by {\cite[Theorem 4.2]{kell3}}, there exists a singular inner function $\theta$ 
such that $\operatorname{ran}\theta(S_{\omega'})$ is not dense. The conclusion of the theorem follows from 
 Lemma~\ref{L:1.1}.

If  condition (i) is fulfilled, then, by {\cite[Theorem 4.1]{kell3}}, there exists a singular inner function $\theta$ 
such that $\operatorname{ran}\theta(S_{w'})$ is not dense. It follows from \eqref{E:7.16} that $S_{\omega'}\prec S_{w'}$. We obtain that 
$T\prec  S_{w'}$.  It remains to apply Lemma~\ref{L:1.1}.  \end{proof}

 \begin{example}\label{E:e7.6} For $a>1$ put $w(-n)=n^{(\log\log n)^a}$ for sufficiently large $n$. Then it is possible to define $w(-n)$ 
for small $n\geq 1$ in such a way that $\{w(n)\}_{n\leq -1}$ satisfies the conditions on $w$ from Theorem~\ref{T:7.5}(i).\end{example}

The following theorem is now a straightforward consequence of {\cite[Theorem 1.1]{bovo}}.

\begin{theorem}\label{T:7.7}
 Suppose that  a sequence $\{p(n)\}_{n=1}^\infty$ of positive numbers has the following properties: $\{p(n)\}_{n=1}^\infty$ is concave,
$\displaystyle\sum_{n\geq 1}\textstyle \frac{p(n)}{n^2}=\infty$, $\displaystyle\lim_{n\to\infty}\textstyle \frac{p(n)}{n} = 0$, there exists $\varepsilon>0$ such that the sequence 
$ \bigl\{\frac{p(n)}{n^\varepsilon}\bigr\}_{n=1}^\infty$ is increasing, and there exists $C>0$ such that 
\begin{align}\label{E:7.17}\log(n+1)\leq Cp(n) \ \text{  for }n\geq 1. \end{align}
Suppose that $\omega$ is a submultiplicative dissymmetric weight, $\omega$ satisfies \eqref{E:7.8}, and 
\begin{align}\label{E:7.18}\lim_{n\to\infty} \frac{\log\omega(-n)}{p(n)} = \infty.  \end{align} 
Let $\tau\subset\mathbb T$ be a Borel set such that $0<m(\tau)<1$, and let $T$ be the operator from Proposition~\ref{P:7.4}. 
Then $T$ is similar to a contraction, $T$ is quasianalytic, $\pi(T)=\tau$, and $\sigma(T)=\mathbb T$.\end{theorem}

\begin{proof} Taking into account Theorem~\ref{T:7.2} and Proposition~\ref{P:7.4}, it remains to prove that if 
$f\oplus u\in P^2(m_{2,\alpha}+\mu)\oplus\ell^2_{\omega -}$ and there exists a measurable set $\tau'\subset\tau$ 
such that $m(\tau')>0$ and \begin{align*}X (f\oplus u)=0 \ \ \ \text{ a.e. on } \tau', \end{align*} 
where $X$ is given in Proposition~\ref{P:7.4}, then $f\oplus u=0$.
Taking into account the definition of $X$,  it is sufficient to verify the following. 
Suppose that  $f\in P^2(m_{2,\alpha})$, $f$ has nontangential boundary values $f(\zeta)$ for a.e. $\zeta\in\tau'$, $g\in H^2_-$, 
$\sum_{n\leq -1}|\widehat g(n)|^2\omega(n)^2 <\infty$, and $f(\zeta)=-g(\zeta)$ for a.e. $\zeta\in\tau'$. 
Then $f\equiv 0$ and $g\equiv 0$.

Let  $f\in P^2(m_{2,\alpha})$. Then $\sum_{n\geq 0}|\widehat f(n)|^2v_\alpha(n)^2 <\infty$ (see \eqref{E:7.12}, \eqref{E:7.13}). Therefore, 
there exists $0<c<1$ such that $|\widehat f(n)|v_\alpha(n)\leq c$ for sufficient large $n$. We have 
$$\log |\widehat f(n)|\leq \log c -\log v_\alpha(n) =\log c +\frac{\alpha+1}{2}\log(n+1)\leq \frac{\alpha+1}{2}Cp(n)$$
for sufficient large $n$ (we apply \eqref{E:7.17} and the estimate $\log c<0$). Thus, 
\begin{align}\label{E:7.20} \limsup_{n\to \infty}\frac{\log|\widehat f(n)|}{p(n)}<\infty. \end{align} 

Let $g\in H^2_-$ be such that  $\sum_{n\leq -1}|\widehat g(n)|^2\omega(n)^2 <\infty$. Again, there exists $0<c<1$ such that 
$|\widehat g(-n)|\omega(-n)\leq c$ for sufficient large $n$. We have 
$$\log |\widehat g(-n)|\leq \log c -\log \omega(-n) \leq -\log \omega(-n).$$ By \eqref{E:7.18},
\begin{align}\label{E:7.21} \lim_{n\to \infty} \frac{\log|\widehat g(-n)|}{p(n)} = - \infty.  \end{align} 

Now we are ready to apply  {\cite[Theorem 1.1]{bovo}}. We have a set $\tau'\subset\mathbb T$ such that  $m(\tau')>0$, a function $f$ analytic in $\mathbb D$ 
such that $f$ satisfies \eqref{E:7.20} and  $f$ has nontangential boundary values $f(\zeta)$ for a.e. $\zeta\in\tau'$, a function $g\in H^2_-$ such that
$g$  satisfies \eqref{E:7.21}, and we have $f(\zeta)=-g(\zeta)$ for a.e. $\zeta\in\tau'$. By  {\cite[Theorem 1.1]{bovo}}, $f\equiv 0$ and $g\equiv 0$. 
\end{proof}

 \begin{example}\label{E:e7.8} For $0<\beta<\beta'\leq 1$ set $\omega(-n)=\exp\bigl(\frac{n}{(\log n+1)^\beta}\bigr)$, $n\geq 1$, $\omega(n)=1$, $n\geq 0$, 
and $p(n)=\frac{n}{(\log n+1)^{\beta'}}$, $n\geq 1$. Then $\omega$ and $\{p(n)\}_{n=1}^\infty$ satisfy the conditions of 
Theorem~\ref{T:7.7}.\end{example}

\subsection*{Acknowlegement.} The author is grateful to L\'aszl\'o K\'erchy for valuable discussions and to an anonymous  referee for pointing out Lemma~\ref{L:2.5} and improving  Corollaries~\ref{C:5.10} and \ref{C:6.10}.

\end{document}